\documentclass[10pt]{amsart}
\usepackage{amssymb,amsmath,latexsym,amscd,amsfonts,amsthm, yfonts}
\usepackage{pgf,tikz,pgfplots}
\usepackage{subcaption}
\pgfplotsset{compat=1.12}
\usetikzlibrary{arrows}
\usetikzlibrary[patterns]
\usepackage{graphics}
\usepackage{epsfig}
\usepackage{psfrag}
\usepackage{wasysym}

\usepackage[none,bottom]{draftcopy}

\def\ignore #1 {}

\newtheorem{thm}{Theorem}
\newtheorem*{thm*}{Theorem}
\newtheorem{lem}[thm]{Lemma}

\newtheorem{defn}[thm]{Definition}
\newtheorem{prop}[thm]{Proposition}
\newtheorem{cor}[thm]{Corollary}
\newtheorem{conjecture}{Conjecture}
\newtheorem{question}{Question}
\newtheorem*{question*}{Question}
\theoremstyle{definition}

\newtheorem{example}{Example}

\newtheorem{questionnuminner}{Question}
\newenvironment{questionnum}[1]{%
  \renewcommand\thequestionnuminner{#1}%
  \questionnuminner
}{\endquestionnuminner}

\def\F{\mbox{{\bf F}}}

\def\Q{\mbox{{\bf Q}}}
\def\R{\mbox{{\bf R}}}
\def\C{\mbox{{\bf C}}}
\newcommand{\CC}{{\mathcal{C}}}
\def\P{\mbox{{\bf P}}}
\def\Z{\mbox{{$\mathbb{Z}$}}}

\def\a{\mbox{${\alpha}$}}

\def\pp{\mathfrak{p}}
\def\qq{\mathfrak{q}}

\def\PP{\mathbf{p}}
\def\QQ{\mathbf{q}}
\def\RR{\mathbf{r}}
\def\SS{\mathbf{s}}

\newcommand{\Diss}{{\mathcal{D}}}
\newcommand{\ftilde}{{\tilde{f}}}

\DeclareMathOperator{\im}{Im}

\DeclareMathOperator{\Aff}{Aff}
\DeclareMathOperator{\Area}{Area}
\DeclareMathOperator{\Vxs}{Vertices}
\DeclareMathOperator{\Triangles}{Triangles}
\DeclareMathOperator{\Constraints}{Constraints}
\DeclareMathOperator{\valence}{Valence}
\DeclareMathOperator{\length}{Length}

\newcommand{\D}{\Delta}
\renewcommand{\a}{\alpha}

\newcommand{\s}{\sigma}

\newcommand{\T}{{\mathcal{T}}}

\def \good {small}
\def \goodness {smallness}

\begin{document}

\title{Generalized dissections and Monsky's Theorem}
 \author{Aaron Abrams
 \and Jamie Pommersheim}

\begin{abstract} 
Monsky's celebrated equidissection theorem follows from his more general proof of the existence of a polynomial relation $f$ among the areas of the triangles in a dissection of the unit square.  More recently, the authors studied a different polynomial $p$, also a relation among the areas of the triangles in such a dissection, that is invariant under certain deformations of the dissection.  In this paper we study the relationship between these two polynomials.

We first generalize the notion of dissection, allowing triangles whose orientation differs from that of the plane. We define a deformation space of these generalized dissections and we show that this space is an irreducible algebraic variety.  We then extend the theorem of Monsky to the context of generalized dissections, showing that Monsky's polynomial $f$ can be chosen to be invariant under deformation. Although $f$ is not uniquely defined, the interplay between $p$ and $f$ then allows us to identify a canonical pair of choices for the polynomial $f$.  In many cases, all of the coefficients of the canonical $f$ polynomials are positive.  We also use the deformation-invariance of $f$ to prove that the polynomial $p$ is congruent modulo 2 to a power of the sum of its variables. 
\end{abstract}

\address{{\tt abramsa@wlu.edu}:  Mathematics Department, Washington and Lee University, Lexington VA 24450, USA}
\address{{\tt jamie@reed.edu}:  Department of Mathematics, Reed College, 3203 SE Woodstock Blvd, Portland OR 97202, USA} 

\keywords{Triangulation, area relation, equidissection}

\maketitle

\section{Introduction}

\bigskip

In 1970 Paul Monsky proved the following theorem:

\begin{thm*}[Monsky \cite{monsky}]
Fix a dissection of the unit square into $n$ triangles, and denote the areas of the triangles by $a_1,\ldots,a_n$.  Then there is an integer polynomial $f$ in $n$ indeterminates such that 
$f(a_1,\ldots,a_n)=1/2.$
\end{thm*}

A corollary is Monsky's famous ``equidissection'' theorem:  if a square is dissected into $n$ triangles of equal area, then $n$ must be even.  This follows because there is no integer polynomial in $n$ variables with $f(\frac 1n,\ldots,\frac 1n)=\frac 12$ when $n$ is odd.

\bigskip

Happy 50th birthday, Monsky's Theorem!

\bigskip

In the half-century since its publication, the equidissection theorem has inspired a significant amount of mathematics, including numerous other equidissection theorems in the plane, higher dimensional analogs, approximation theorems, and more.  Relatively little attention has been focused on the polynomial $f$, however.

A \emph{dissection} of a square is defined as a finite collection of triangles in the plane whose interiors do not intersect and whose union is the square.  Monsky's theorem is a statement about dissections.  

Over the years it has occurred to several people to first fix the combinatorics of a dissection, and then try to understand which collections of areas \emph{are} realized by the triangles.  We heard of this approach from Joe Buhler, whose student Adam Robins wrote \cite{transcendence} about it, and from Serge Tabachnikov, whose students Joshua Kantor and Max Maydanskiy wrote \cite{gonewild} about it.
In this direction, we established in \cite{triangles1} the existence of a nonzero integer polynomial $p$, different than $f$, associated to certain dissections which also has one variable for each triangle and which vanishes, rather than taking the value $1/2$, on the input $(a_1,\ldots,a_n)$.
By construction $p$ depends only on the combinatorics of the dissection, so the same $p$ also vanishes at any tuple of areas arising by deforming the dissection; indeed under the hypotheses of our theorem the zero set of $p$ is exactly the \emph{area variety} of the triangulation, which is the (closure of the) collection of realizable areas.
For elementary reasons $p$ is irreducible and homogeneous.

The polynomials $p$ and $f$ are our primary objects of study.  These polynomials are closely related, although they have different roles in the theory.  In \cite{triangles1} we called $p$ a \emph{Monsky polynomial}, but here we emphasize the distinction and the interplay between the two, so we give them different names:  $p$ is the \emph{area polynomial} and $f$ is the \emph{Monsky polynomial}.

Here are two examples which make numerous appearances throughout the paper.
\begin{example}\label{ex:1}
The dissection in Figure \ref{fig:intro1} has area polynomial $p=A-B+C-D$ (or its negative), and Monsky polynomial $f=A+C$ (or $\tilde f=B+D$, or $f+p=2A-B+2C-D$, etc.).  One can easily see that regardless of where the central vertex is placed, the polynomial $p$ evaluates to zero, and as long as the square has unit area, $f, \tilde f,$ and $f+p$ all evaluate to $1/2$.
For any square, $f$ evaluates to half the total area.
\end{example}
\begin{example}\label{ex:2}
Less apparently, the dissection in Figure \ref{fig:intro2} has 
\begin{align*}
p &= A^2+C^2+E^2-2AC+2AE+2CE-B^2-D^2-F^2-2BD-2BF+2DF \\
\mbox{and } f &= A^2+C^2+E^2+2AE+2CE+2DF +(A+C+E)(B+D+F).
\end{align*}
Again, regardless of the placement of $u$ and $v$, $p$ evaluates to zero and $f$ evaluates to half the area of the square.
\end{example}

\begin{figure}
\centering
    \begin{subfigure}[t]{.4\linewidth}
    \subcaptionbox{Example 1: \\ $p=A-B+C-D$ \\ $f=A+C$\label{fig:intro1}}
    {
        \begin{tikzpicture}[scale=.5]
        \draw (0,0) -- (0,6) -- (6,6) -- (6,0) -- cycle;
\draw (4,2.5) -- (0,0);
\draw (4,2.5) -- (6,0);
\draw (4,2.5) -- (0,6);
\draw (4,2.5) -- (6,6);
\draw [fill] (0,0) circle (1pt);
\draw [fill] (0,6) circle (1pt);
\draw [fill] (6,6) circle (1pt);
\draw [fill] (6,0) circle (1pt);
\draw [fill] (4,2.5) circle (1pt);
\draw (3.5,1) node{$A$};
\draw (5,2.8) node{$B$};
\draw (3.5,4.5) node{$C$};
\draw (1.5,2.5) node{$D$};
        \end{tikzpicture}   
        }
    \end{subfigure}
\qquad
    \begin{subfigure}[t]{.4\linewidth}
    \subcaptionbox{Example 2: \\ See text for $p$ and $f$. \label{fig:intro2}}
    {
        \begin{tikzpicture}[scale=.5]
        \draw (0,0) -- (0,6) -- (6,6) -- (6,0) -- cycle;
\draw (0,0) -- (2.4,1.6) -- (4.4,3.6) -- (6,6);
\draw (6,0) -- (2.4,1.6) -- (0,6) -- (4.4,3.6) -- cycle;
\draw[fill=black] (0,0) circle (1pt);
\draw[fill=black] (0,6) circle (1pt);
\draw[fill=black] (6,0) circle (1pt);
\draw[fill=black] (6,6) circle (1pt);
\draw[fill=black] (2.4,1.6) circle (1pt);
\draw[fill=black] (4.4,3.6) circle (1pt);
\draw (3,.5) node{$A$};
\draw (4.2,1.8) node{$B$};
\draw (5.4,2.8) node{$C$};
\draw (1,2.2) node{$D$};
\draw (2.6,3.4) node{$E$};
\draw (3.8,5) node{$F$};
\draw (2.4,1.6) node[anchor=north]{$u$};
\draw (4.4,3.6) node[anchor=west]{$v$};
        \end{tikzpicture}   
        }
    \end{subfigure}
    \caption{Seminal examples:  $p$ evaluates to zero and $f$ evaluates to $1/2$.}
    \label{fig:intro}
\end{figure}

The fact that $p$ is well-defined\footnote{up to sign} essentially reflects the correctness of a heuristic dimension count, whereas Monsky's polynomial $f$ provides number-theoretic (specifically, mod 2) information.
However, in Monsky's theorem a dissection is treated as a static object, and invariance of $f$ under deformation is not guaranteed.  
Think of a dissection in which some triangles in the middle, say $i$ and $j$, have areas summing to $1/2$. 
Then the polynomial $f(x_1,\ldots,x_n)=x_i+x_j$ satisfies the conclusion of Monsky's theorem, but the sum $a_i+a_j$ could easily change when the dissection is deformed.  One of our main goals is to extend Monsky's theorem to show that $f$ can be made deformation-invariant, as it is in the examples we have already seen.

\begin{figure} 
\begin{center} 
\begin{tikzpicture}[scale=.5]
    \draw (0,0) -- (0,6) -- (6,6) -- (6,0) -- cycle;
\draw (0,2) -- (2,6) -- (6,4) -- (4,0) -- cycle;
\draw (0,2) -- (6,4) ;
\draw[fill=black] (0,0) circle (1pt);
\draw[fill=black] (0,6) circle (1pt);
\draw[fill=black] (6,0) circle (1pt);
\draw[fill=black] (6,6) circle (1pt);
\draw[fill=black] (0,2) circle (1pt);
\draw[fill=black] (2,6) circle (1pt);
\draw[fill=black] (6,4) circle (1pt);
\draw[fill=black] (4,0) circle (1pt);
\end{tikzpicture} 
\qquad
        \begin{tikzpicture}[scale=.5]
        \draw (-4,6) -- (6,6) -- (6,0) -- (0,0) -- (0,6);
\draw (-4,6) -- (6,0);
\draw[fill=black] (-4,6) circle (1pt);
\draw[fill=black] (6,6) circle (1pt);
\draw[fill=black] (6,0) circle (1pt);
\draw[fill=black] (0,0) circle (1pt);
\draw[fill=black] (0,6) circle (1pt);
\draw[fill=black] (0,3.6) circle (1pt);
        \end{tikzpicture}   
\caption{A dissection (left) and a generalized dissection (right) of a square.}
\label{fig:samples}
\end{center} 
\end{figure}
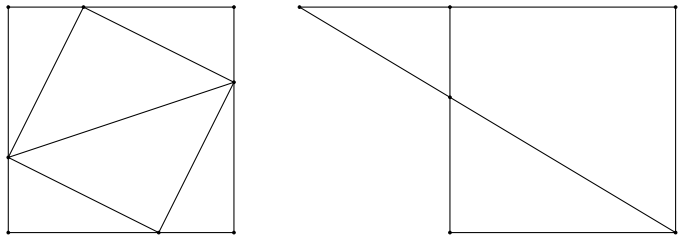

It turns out that the act of deforming a dissection is trickier than it may appear at first glance, and it deserves to be taken seriously.  
One issue is that the vertices may be constrained to lie on certain line segments, so in general the vertices cannot move freely and independently of each other.  Another issue is that one is forced to confront the possibility that triangles might degenerate, or turn upside-down. 
In the first part of this paper we develop a framework for handling these issues, building on our work in \cite{triangles1}.  
The main idea is to view a dissection as the image of a certain map which itself has a natural deformation space.
We are led to a notion of a \emph{generalized dissection}, and we will see that there are generalized dissections that cannot be deformed back into (classical) dissections. See Figure \ref{fig:samples} (right), where there are three triangles, one of which is upside-down.  Examples like this turn out to be crucial to our theory. 

Our main theorem about these deformation spaces, which we call $X$, is that they are irreducible rational varieties.  This is proved in Section \ref{sec:X}.
Our proof is more subtle than we anticipated, because we encountered fundamental issues about arrangements of points and lines that required some finesse to mitigate.  In Section \ref{sec:rambling} we discuss some questions that arose in this process and their relationship with the well-studied areas of point/line configurations and oriented matroids.

A different method for treating the problem of deformations has been proposed and studied in ~\cite{ziegler} by Labb\'e, Rote, and Ziegler, who were interested in approximating equidissections.

In the second part of the paper, with the foundations now established, we are able to investigate $p$ and $f$.  Our main technical results (Theorems Monsky+ and Monsky++) extend Monsky's theorem to the deformation space $X$, showing that $f$ can indeed be chosen to be invariant under deformation.
That is, we give algebraic versions of the theorem, showing that for any (generalized) dissection, not only do the \emph{areas} of the triangles satisfy a polynomial relation, but also the \emph{formulas} for the areas satisfy a polynomial relation.  Thus we may think of $f$ as a ``dynamic'' object, as we did already with $p$.

Once both $p$ and $f$ are thusly defined, we are finally able to rigorously explore the relationships between the two.  In Sections \ref{sec:mod2} and \ref{sec:canonical} we prove our main results about $p$ and $f$ by exploiting features of each polynomial to deduce information about the other.

Specifically, recall from Example 1 above that a given dissection has many Monsky polynomials.  The canonicalness of $\{p,-p\}$ allows us to define a canonical pair $\{f,\tilde f\}$ with extra known and conjectured properties; for example these have minimal degree among deformation-invariant polynomials satisfying Monsky's theorem.  
These polynomials also often have non-negative coefficients, an observation we will return to in Section \ref{sec:positivity}. 

In the other direction the number-theoretic content of $f$ transports to $p$, giving additional information about its structure.
For instance we show that mod 2, the polynomial $p$ is congruent to a power of the sum of the variables.

We close in Section \ref{sec:equi} with a question about equidissections.

One of the pleasant features of the present setup is that we minimize the amount of combinatorial information needed to parameterize the deformation space of a generalized dissection.  This information is often implicit in a drawing of the dissection, and this setup simplifies the computation of $p$ and $f$ relative to what we did in \cite{triangles1}. 
The job is still inherently computationally expensive, but the cost now essentially depends only on the number of triangles in the dissection, and not how much degeneracy there is.

Many mysteries remain about these polynomials.

\tableofcontents

\part{Deforming dissections}

\bigskip

In this part of the paper we develop the language of generalized dissections and constrained triangulations, which we use to define deformation spaces of dissections.

\section{Generalized dissections}

Classically, a \emph{dissection} of a square is a finite collection of triangles in the (Euclidean) plane whose interiors do not intersect and whose union is the square.  Our first goal here is to give a more general definition that allows for deformations.  We start by setting some terminology.

We work in the affine plane $\C^2$.  (The reader who prefers to think of everything taking place in $\R^2$ is encouraged to do so; we prove in Section \ref{sec:really} that this makes no difference to our theory.)

If $S$ is a cyclically ordered finite set $S=(s_1,\ldots,s_n)$ then we define an \emph{edge} of $S$ to be any of the ordered pairs $(s_i,s_{i+1})$, with indices taken mod $n$.

A \emph{polygon}, or \emph{$n$-gon}, is a cyclically ordered set of $n\ge3$ distinct points in $\C^2$, called vertices.  A $3$-gon is also called a \emph{triangle}; thus a triangle comes with an orientation.  
A polygon is \emph{totally degenerate} if its vertices are collinear, \emph{degenerate} if it has three consecutive vertices (in the cyclic order) that are collinear, and \emph{non-degenerate} if no three consecutive vertices are collinear.

An \emph{abstract polygon}, or \emph{abstract $n$-gon}, is a $2$-cell whose boundary circle consists of $n$ $0$-cells (vertices) and $n$ $1$-cells (also called edges).

Corresponding to any polygon $\pentagon$ (including degenerate and totally degenerate ones) is an abstract polygon whose vertices are labeled by the points of $\pentagon$ (in the same cyclic order).
Associated to a family of polygons we can construct an abstract $2$-dimensional complex from a corresponding family of abstract polygons by gluing together along edges: the edge $(v,w)$ of one polygon is glued to the edge $(w,v)$ of another.

Notice that we have chosen to label the vertices of the abstract polygons and complexes by the points themselves.  For example, the vertex of the abstract polygon corresponding to $(1,1)$ is called $(1,1)$.\footnote{Another reasonable name for this vertex would have been $v_{(1,1)}$, which has the advantage of emphasizing the abstract nature of this vertex, but the disadvantage of being clearer.} 

\begin{defn}\label{def:generalized}
Let $\pentagon$ be a polygon in $\C^2$.   A \emph{generalized dissection} of $\pentagon$ consists of a finite set $\Triangles$ and a finite set $\Constraints$ such that:
\begin{enumerate}
    \item Each element of $\Triangles$ is a non-degenerate triangle in $\C^2$.
    \item Each element of $\Constraints$ is a totally degenerate polygon in $\C^2$, each of whose vertices is a vertex of at least one triangle in $\Triangles$
    \item Any two distinct constraints share at most one vertex
    \item The associated 2-complex built from abstract polygons corresponding to the union of $\Triangles$ and $\Constraints$ is an oriented disk with boundary equal to $\pentagon$.
\end{enumerate}
\end{defn}

Note that $\pentagon$ is allowed to be degenerate or totally degenerate. 

We think of the elements of $\Triangles$ as the triangles in the dissection, except now they are oriented.  
Elements of $\Constraints$ are interpreted as collinearity constraints; item (3) ensures that the constraints are maximal.  
The abstract polygons corresponding to elements of $\Constraints$ are called \emph{poofagons}.
(These may or may not be triangles, but they are not elements of $\Triangles$.)
Item (4) implies that we can interpret the data as the image of a PL map from a cellulated disk into the plane, under which the poofagons have degenerated into line segments.
(Not every such map gives a generalized dissection though, as illustrated below by Figure \ref{fig:non-dissection}.)

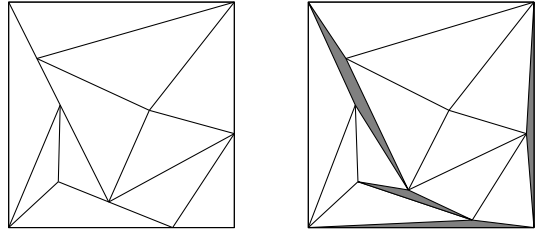
\begin{figure} 
    \begin{subfigure}[t]{.25\linewidth}
        \begin{tikzpicture}[scale=.5, rotate=90]
\draw (0,0) -- (6,0) -- (6,6) -- (0,6) -- cycle;
\draw  (0,0)-- (6,0);
\draw  (6,0)-- (6,6);
\draw (6,6)-- (0,6);
\draw (0,6)-- (0,0);
\draw (2.50,0)-- (0,1.64);
\draw  (0,1.64)-- (1.22,4.68);
\draw  (1.22,4.68)-- (0,6);
\draw  (0.68,3.34)-- (6,6);
\draw  (1.22,4.68)-- (3.26,4.63);
\draw  (3.26,4.63)-- (0,6);
\draw  (4.50,5.25)-- (6,0);
\draw  (0.68,3.34)-- (2.50,0);
\draw  (2.50,0)-- (3.12,2.26);
\draw  (3.12,2.26)-- (0.68,3.34);
\draw  (4.50,5.25)-- (3.12,2.26);
\draw  (3.12,2.26)-- (6,0);
        \end{tikzpicture} 
    \end{subfigure}
\qquad
    \begin{subfigure}[t]{.25\linewidth}
        \begin{tikzpicture}[scale=.5, rotate=90]
\draw (0,0) -- (6,0) -- (6,6) -- (0,6) -- cycle;
\draw  (0,0)-- (6,0);
\draw  (6,0)-- (6,6);
\draw (6,6)-- (0,6);
\draw (0,6)-- (0,0);
    \draw (2.50,.2) -- (.2,1.64);
    \draw [fill=gray] (2.50,.2) -- (0,0) -- (6,0) -- cycle;

    \draw  (.2,1.64)-- (1.22,4.68);
    \draw [fill=gray] (.2,1.64)-- (0,6) -- (0,0) -- cycle;

\draw  (1.22,4.68)-- (0,6);
    \draw (1,3.34)-- (6,6);
    \draw [fill=gray] (1,3.34) -- (.2,1.64) -- (1.22,4.68) -- cycle;
    
\draw  (1.22,4.68)-- (3.26,4.75);
    \draw [fill=gray] (3.26,4.75) -- (6,6) -- (4.50,5) -- (1,3.34)-- cycle;
\draw  (3.26,4.75)-- (0,6);
\draw  (4.50,5)-- (6,0);
\draw  (1,3.34)-- (2.50,.2);
\draw  (2.50,.2)-- (3.12,2.26);
\draw  (3.12,2.26)-- (1,3.34);
\draw  (4.50,5)-- (3.12,2.26);
\draw  (3.12,2.26)-- (6,0);
        \end{tikzpicture} 
    \end{subfigure}
\caption{A dissection and its associated 2-complex.  The shaded cells are the poofagons.}
\label{fig:poofing}
\end{figure}

Often, the sets $\Triangles$ and $\Constraints$ are implicitly defined by a drawing.  For classical dissections this is always the case, as we prove in Proposition \ref{prop:DtoDiss} below.  An example of a classical dissection is shown in Figure \ref{fig:poofing}, along with the 2-complex associated to the corresponding (implicitly defined) generalized dissection.   The generalized dissection has four poofagons, one quadrilateral and the rest triangles, shown as shaded cells.

\begin{figure}
    \begin{subfigure}[t]{.4\linewidth}
    \subcaptionbox{\label{fig:ACE-D}}
    {
        \begin{tikzpicture}[scale=.4]
        \draw (-4,6) -- (6,6) -- (6,0) -- (0,0) -- (0,6);
\draw (-4,6) -- (6,0);
\draw[fill=black] (-4,6) circle (1pt);
\draw[fill=black] (6,6) circle (1pt);
\draw[fill=black] (6,0) circle (1pt);
\draw[fill=black] (0,0) circle (1pt);
\draw[fill=black] (0,6) circle (1pt);
\draw[fill=black] (0,3.6) circle (1pt);
        \end{tikzpicture}   
        }
    \end{subfigure}
    \begin{subfigure}[t]{.4\linewidth}
    \subcaptionbox{\label{fig:ACE-poofed}}
    {
        \begin{tikzpicture}[scale=.4]
        \draw (0,0) -- (0,6) -- (6,6) -- (6,0) -- cycle;
\draw (0,0) -- (2.4,1.6) -- (4.4,3.6) -- (6,6);
\draw (6,0) -- (2.4,1.6) -- (0,6) -- (4.4,3.6) -- cycle;
\draw[fill=black] (0,0) circle (1pt);
\draw[fill=black] (0,6) circle (1pt);
\draw[fill=black] (6,0) circle (1pt);
\draw[fill=black] (6,6) circle (1pt);
\draw[fill=black] (2.4,1.6) circle (1pt);
\draw[fill=black] (4.4,3.6) circle (1pt);
\draw[fill=gray] (0,0) -- (0,6) -- (2.4,1.6) -- cycle;
\draw[fill=gray] (0,6) -- (4.4,3.6) -- (6,6) -- cycle;
\draw[fill=gray] (6,0) -- (4.4,3.6) -- (2.4,1.6) -- cycle;
        \end{tikzpicture}   
        }
    \end{subfigure}
\caption{The $ACE$ example:  a generalized dissection with constraints.  This cannot be deformed to a dissection.}
    \label{fig:ACE}
\end{figure}
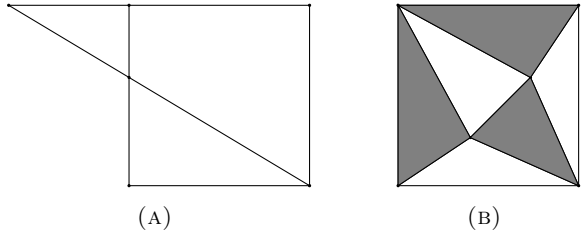

One should acquaint oneself with a few more examples before proceeding.  Some basic ones are shown in Figure \ref{fig:examples-D}, and we separately highlight an especially important one in Figure \ref{fig:ACE}.

\begin{example}[cf.~Example \ref{ex:1}]\label{ex:1-D}
Figure \ref{fig:diag1} is a dissection with four triangles; it is also a generalized dissection with the same four triangles (now oriented) and no constraints.  The corresponding abstract triangles glue together to form a simplicial complex of which this is a drawing.

Figure \ref{fig:diag1-gen} resembles \ref{fig:diag1}, except the central vertex has been dragged outside the square.  Here and elsewhere, we have indicated the vertices with small dots in order to avoid potential confusion with edges that intersect at points of the plane that are not vertices.  This is not a dissection.  It is a generalized dissection with four triangles, one of which is oriented differently from the plane.  As in Figure \ref{fig:diag1}, there are no constraints.  The corresponding abstract triangles form the same simplicial complex as Figure \ref{fig:diag1}.
\end{example}

\begin{example}[cf.~Example \ref{ex:2}]\label{ex:2BE}
Figure \ref{fig:BE-D} is a dissection, both classical and generalized, with four triangles. The generalized dissection has a constraint, which is a (totally degenerate) quadrilateral.  The associated cell complex can be triangulated in two ways by choosing a diagonal of this quadrilateral; one of the resulting simplicial complexes is shown later in Figure \ref{fig:diag2}.  
\end{example}

\begin{example}
Figure \ref{fig:non-dissection} is not a generalized dissection at all.  Although its faces ``cancel,'' the flattened tetrahedron pinned to the center of the square makes it impossible to describe this as a generalized dissection.  In particular, the simplicial complex made from the obvious eight (abstract) triangles is homeomorphic to the one-point union of a disk and a sphere.
\end{example}

\begin{example}[The $ACE$ example, cf.~Example \ref{ex:2}]
\label{ex:2ACE-D}
Finally, Figure \ref{fig:ACE-D} is a generalized dissection with three triangles and three constraints.  It is an interesting specimen.  It takes a moment to identify the triangles and the constraints (with the correct orientations).  There are three triangles, one of which is upside-down.  The reader should verify that this this does indeed satisfy the definitions of a generalized dissection, with the associated simplicial complex shown in Figure \ref{fig:ACE-poofed}, with the poofagons shaded.  (This is the same 2-complex shown later, in Figure \ref{fig:diag2}.)  One feature of this example is that it cannot be deformed into a (classical) dissection in which all three triangles remain alive.
\end{example}

\begin{figure}
    \begin{subfigure}[t]{.25\linewidth}
    \subcaptionbox{A dissection. \label{fig:diag1}}
    {
        \begin{tikzpicture}[scale=.3]
        \draw (0,0) -- (0,6) -- (6,6) -- (6,0) -- cycle;
\draw (4,2.5) -- (0,0);
\draw (4,2.5) -- (6,0);
\draw (4,2.5) -- (0,6);
\draw (4,2.5) -- (6,6);
\draw [fill] (0,0) circle (1pt);
\draw [fill] (0,6) circle (1pt);
\draw [fill] (6,6) circle (1pt);
\draw [fill] (6,0) circle (1pt);
\draw [fill] (4,2.5) circle (1pt);
        \end{tikzpicture}   
        }
    \end{subfigure}
\qquad
    \begin{subfigure}[t]{.25\linewidth}
    \subcaptionbox{A generalized dissection that is not a dissection. \label{fig:diag1-gen}}
    {
        \begin{tikzpicture}[scale=.3]
        \draw (0,0) -- (0,6) -- (6,6) -- (6,0) -- cycle;
\draw (7.5,2.5) -- (0,0);
\draw (7.5,2.5) -- (6,0);
\draw (7.5,2.5) -- (0,6);
\draw (7.5,2.5) -- (6,6);
\draw [fill] (0,0) circle (1pt);
\draw [fill] (6,0) circle (1pt);
\draw [fill] (0,6) circle (1pt);
\draw [fill] (6,6) circle (1pt);
\draw [fill] (7.5,2.5) circle (1pt);
        \end{tikzpicture}   
        }
    \end{subfigure}
\\[10pt]
    \begin{subfigure}[t]{.25\linewidth}
    \subcaptionbox{A dissection with a constraint. \label{fig:BE-D}}
    {
        \begin{tikzpicture}[scale=.3]
        \draw (0,0) -- (0,6) -- (6,6) -- (6,0) -- cycle;
\draw (6,0) -- (0,6) ;
\draw (0,0) -- (2,4);
\draw (4,2) -- (6,6);
\draw [fill] (0,0) circle (1pt);
\draw [fill] (0,6) circle (1pt);
\draw [fill] (6,6) circle (1pt);
\draw [fill] (6,0) circle (1pt);
\draw [fill] (4,2) circle (1pt);
\draw [fill] (2,4) circle (1pt);
        \end{tikzpicture}   
        }
    \end{subfigure}
\qquad
    \begin{subfigure}[t]{.25\linewidth}
    \subcaptionbox{Not a dissection. \label{fig:non-dissection}}
    {
        \begin{tikzpicture}[scale=.3]
\draw (0,0) -- (0,6) -- (6,6) -- (6,0) -- cycle;
\draw (4,2.5) -- (0,0);
\draw (4,2.5) -- (6,0);
\draw (4,2.5) -- (0,6);
\draw (4,2.5) -- (6,6);
\draw (4,2.5) -- (8,3.5) -- (8.8,4.7) -- (7.8,5.7) -- cycle;
\draw (8,3.5) -- (7.8,5.7) ;
\draw[color=gray] (4,2.5) -- (8.8,4.7);
\draw [fill] (0,0) circle (1pt);
\draw [fill] (0,6) circle (1pt);
\draw [fill] (6,0) circle (1pt);
\draw [fill] (6,6) circle (1pt);
\draw [fill] (4,2.5) circle (1pt);
\draw [fill] (8,3.5) circle (1pt);
\draw [fill] (8.8,4.7) circle (1pt);
\draw [fill] (7.8,5.7) circle (1pt);

\draw [white, very thick] (5, 2.75) -- (7, 3.25);
\draw (5, 2.75) -- (7, 3.25);
\draw [white, very thick] (5.2, 3.05) -- (7.6, 4.15);
\draw [gray] (5.2, 3.05) -- (7.6, 4.15);
\draw [white, very thick] (4.95, 3.3) -- (6.85, 4.9);
\draw (4.95, 3.3) -- (6.85, 4.9);
        \end{tikzpicture}   
        }
    \end{subfigure}
\caption{Some basic examples}
\label{fig:examples-D}
\end{figure}
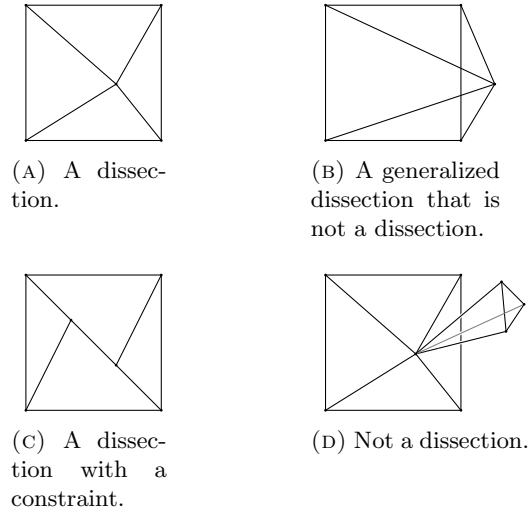

\begin{prop}[$D\leadsto\Diss$]\label{prop:DtoDiss}
The triangles of any (classical) dissection of a square, when  oriented counterclockwise in $\R^2$, comprise the set $\Triangles$ of a generalized dissection.
\end{prop}

\begin{proof}
Let $D$ be a dissection, and let $\Triangles$ be the set of triangles, each oriented counterclockwise.
We need only to specify the collinearity constraints.

Recall that $D$ consists of triangles in $\R^2$.
Say that a vertex of $D$ is \emph{constrained} if it is in the interior of an edge of (a triangle of) $D$ and \emph{unconstrained} otherwise.
Define a \emph{segment} of $D$ to be a line segment in the plane that is contained in the union of the boundaries (edges) of the triangles of $D$ and contains no unconstrained vertex in its interior.
Finally a segment is called \emph{maximal} if it is not contained in any larger segment and it contains at least three vertices of $D$. 

For each maximal segment $M$, we define a constraint containing exactly those vertices that are contained in $M$.
To determine the cyclic order, we use the fact that every  vertex $v$ in the interior of $M$ is constrained, so all edges containing such $v$ (and not contained in $M$) are on the same side of $M$.
Precisely, we traverse the boundary of a small regular neighborhood of $M$ in the plane, counterclockwise.  
Each time we cross an edge of $D$, we record the vertex in $M$ that the edge contains.  
After eliminating duplicates, we have a cyclic ordering on the vertices contained in $M$.
The set $\Constraints$ consists of the cyclically ordered sets constructed in this way.

It is now easy to see that we have a generalized dissection.  
Item (3) of the definition is satisfied since constraints intersect exactly where the corresponding maximal segments intersect, and 
two such segments cannot overlap in an interval by maximality.
Item (4) of the definition is also satisfied because the associated $2$-complex is made by cutting the square open along the maximal segments and gluing in poofagons corresponding to the constraints.
\end{proof}

\begin{defn}
A generalized dissection is \emph{generic} if no line in $\C^2$ contains two intersecting constraints.
\end{defn}

\begin{example}
Figure \ref{fig:non-generic} shows a non-generic dissection $D$ on the left.
The vertex $v$ is unconstrained; our definition of generalized dissection does not allow us to interpret the entire horizontal segment containing $v$ as a single constraint.  Instead we view this segment as two separate constraints intersecting at $v$; this violates the definition of generic.  The middle and right figures show two generic dissections that are close to $D$.  The middle figure has two constraints, whereas on the right the constraints have been merged and there is an additional vertex.
We will see in Section \ref{sec:X} that these two generic dissections have different deformation spaces.

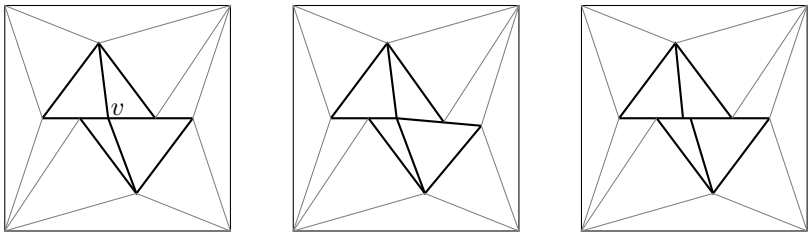
\begin{figure} 
\begin{center}
    \begin{tikzpicture}[scale=.5]
    \draw (0,0) -- (0,6) -- (6,6) -- (6,0) -- cycle;

\draw[thick] (1,3) -- (2.5,5) -- (4,3) ;
\draw[thick] (2,3) -- (3.5,1) -- (5,3) ;
\draw[thick] (1,3) -- (5,3) ;
\draw[thick] (2.5,5) -- (2.75,3) -- (3.5,1) ;

\draw (3,3.25) node{$v$};

\draw[gray] (0,0) -- (1,3) -- (0,6) -- (2.5,5) -- (6,6) -- (4,3) ;
\draw[gray] (6,6) -- (5,3) -- (6,0) -- (3.5,1) -- (0,0) -- (2,3) ;
    \end{tikzpicture}
    \qquad
    \begin{tikzpicture}[scale=.5]
    \draw (0,0) -- (0,6) -- (6,6) -- (6,0) -- cycle;

\draw[thick] (1,3) -- (2.5,5) -- (4,2.9) ;
\draw[thick] (2,3) -- (3.5,1) -- (5,2.8) ;
\draw[thick] (1,3) -- (2.75,3) -- (5,2.8) ;
\draw[thick] (2.5,5) -- (2.75,3) -- (3.5,1) ;


\draw[gray] (0,0) -- (1,3) -- (0,6) -- (2.5,5) -- (6,6) -- (4,2.9) ;
\draw[gray] (6,6) -- (5,2.8) -- (6,0) -- (3.5,1) -- (0,0) -- (2,3) ;
    \end{tikzpicture}
    \qquad
    \begin{tikzpicture}[scale=.5]
    \draw (0,0) -- (0,6) -- (6,6) -- (6,0) -- cycle;

\draw[thick] (1,3) -- (2.5,5) -- (4,3) ;
\draw[thick] (2,3) -- (3.5,1) -- (5,3) ;
\draw[thick] (1,3) -- (5,3) ;
\draw[thick] (2.5,5) -- (2.7,3) ;
\draw[thick] (2.9,3) -- (3.5,1) ;


\draw[gray] (0,0) -- (1,3) -- (0,6) -- (2.5,5) -- (6,6) -- (4,3) ;
\draw[gray] (6,6) -- (5,3) -- (6,0) -- (3.5,1) -- (0,0) -- (2,3) ;
    \end{tikzpicture}
\caption{A non-generic dissection, and two generic dissections.}
\label{fig:non-generic}
\end{center} 
\end{figure}

\end{example}

In Section \ref{sec:X} we will define a generic drawing, and we will see that the two uses of the term ``generic'' line up.  

\begin{question}\label{q:non-generic}
Is every generalized dissection close to a generic one?
\end{question}

\begin{question}\label{q:non-generic-classical}
Is every dissection close to a generic one?
\end{question}

These are questions of incidence geometry, and the answers may depend on the underlying field.  The meaning of ``close'' will be made precise in Section \ref{sec:X}.

\section{Constrained triangulations}

A generalized dissection of a square has an associated 2-complex which is homeomorphic to a disk.  
If we triangulate any non-triangular poofagons, the result leads to what we call a \emph{constrained\footnote{In our previous paper \cite{triangles1} we referred to this as a \emph{generalized triangulation.}} triangulation}.

\begin{defn} \label{def:constrained} 
A \emph{constrained triangulation} $\T$ is a pair $\T=(T,\CC)$, where $T$ is an oriented simplicial complex homeomorphic to a disk, and where $\CC=\{C_i\}$ is a (finite) set of \emph{(collinearity) constraints}.  
Vertices on the boundary of $T$ are called \emph{corners}, and other vertices of $T$ are called \emph{interior vertices.}  
Each collinearity constraint $C_i$ is a set of vertices of $T$ of the form $\Vxs(S_i)$ where $S_i$ is a contiguous set of triangles of $T$.  
(This means that there is a connected subgraph of the dual graph to $T$ whose vertices are the triangles of $S_i$.)  
We require the sets $S_i$ of triangles to be disjoint, although the constraints $C_i$ need not be.  

A 2-cell of $T$ is called \emph{alive} or \emph{living} if there is no constraint containing
all of its vertices.  

Except in Section \ref{sec:rambling}, a constrained triangulation always has four corners, which are labeled $\PP,\QQ,\RR,\SS$ in the cyclic order determined by the orientation of $T$.
\end{defn}

A note about our usage:  much of the modern and classical literature uses the word ``triangulation'' to distinguish a special type of dissection, namely a simplicial one.  However our usage is different.  We use the word ``dissection,'' modified in various ways, to refer to a concrete (geometric) object, whereas a ``triangulation'' is an abstract (topological) object.  It is helpful to think of a dissection as a \emph{drawing} of a triangulation; in fact we make this precise in Section \ref{sec:X}.  (This is how we will deform a dissection.)  So indeed triangulations are always simplicial but a simplicial dissection, which consists of actual triangles in the plane, is not the same thing as a triangulation, which is an abstract simplicial complex.

\begin{prop}[$\Diss\leadsto\T$]\label{prop:DtoT}
Let $\Diss$ be a generalized dissection. There is a constrained triangulation $\T(\Diss)$ whose vertices and living triangles are in 1-1 correspondence with the vertices and triangles of $\Diss$.
\end{prop}

\begin{proof}
Triangulate the poofagons of the associated 2-complex arbitrarily and for each poofagon define a constraint consisting of the vertices of the poofagon, using the boundary to determine the cyclic order.
\end{proof}

Triangulating the poofagons in a different way produces a (slightly) different $\T$ satisfying the conclusion of the proposition, and any two such $\T$'s are related in this way.

If $\CC$ is empty then $\T$ is an abstract version of a classical simplicial dissection of a square.  We call this an \emph{honest triangulation}.  All triangles in an honest triangulation are alive. The constrained triangulation associated to any classical dissection with no constrained vertices (i.e., what is classically called a ``simplicial dissection'') is honest. 

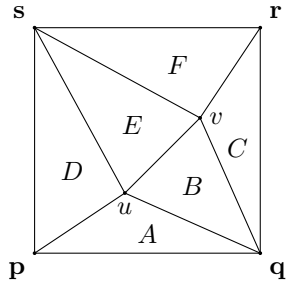
\begin{figure}
    \centering
    \begin{tikzpicture}[scale=.5]
        \draw (0,0) -- (0,6) -- (6,6) -- (6,0) -- cycle;
\draw (0,0) -- (2.4,1.6) -- (4.4,3.6) -- (6,6);
\draw (6,0) -- (2.4,1.6) -- (0,6) -- (4.4,3.6) -- cycle;
\draw[fill=black] (0,0) circle (1pt);
\draw[fill=black] (0,6) circle (1pt);
\draw[fill=black] (6,0) circle (1pt);
\draw[fill=black] (6,6) circle (1pt);
\draw[fill=black] (2.4,1.6) circle (1pt);
\draw[fill=black] (4.4,3.6) circle (1pt);
\draw (3,.5) node{$A$};
\draw (4.2,1.8) node{$B$};
\draw (5.4,2.8) node{$C$};
\draw (1,2.2) node{$D$};
\draw (2.6,3.4) node{$E$};
\draw (3.8,5) node{$F$};
\draw (0,0) node[anchor=north east]{$\PP$};
\draw (6,0) node[anchor=north west]{$\QQ$};
\draw (6,6) node[anchor=south west]{$\RR$};
\draw (0,6) node[anchor=south east]{$\SS$};
\draw (2.4,1.6) node[anchor=north]{$u$};
\draw (4.4,3.6) node[anchor=west]{$v$};
    \end{tikzpicture} 
    \caption{A (drawing of)\protect\footnotemark\ an honest triangulation $T$, with everything labeled.  This is Example 2.}
    \label{fig:diag2}
\end{figure}

\begin{figure}
    \begin{subfigure}[t]{.25\linewidth}
    \subcaptionbox{Two constraints again. \label{fig:AF-T}}
    {
        \begin{tikzpicture}[scale=.4]
        \draw (0,0) -- (0,6) -- (6,6) -- (6,0) -- cycle;
\draw (0,0) -- (2.4,1.6) -- (4.4,3.6) -- (6,6);
\draw (6,0) -- (2.4,1.6) -- (0,6) -- (4.4,3.6) -- cycle;
\draw[fill=black] (0,0) circle (1pt);
\draw[fill=black] (0,6) circle (1pt);
\draw[fill=black] (6,0) circle (1pt);
\draw[fill=black] (6,6) circle (1pt);
\draw[fill=black] (2.4,1.6) circle (1pt);
\draw[fill=black] (4.4,3.6) circle (1pt);
 \draw[fill] (3,.5) circle (2pt);
 \draw[fill] (3.8,5) circle (2pt);
        \end{tikzpicture} 
        }
    \end{subfigure}
\qquad
    \begin{subfigure}[t]{.25\linewidth}
    \subcaptionbox{Here $\CC$ contains one constraint, $\QQ v\SS u$. \label{fig:BE-T}}
    {
        \begin{tikzpicture}[scale=.4]
        \draw (0,0) -- (0,6) -- (6,6) -- (6,0) -- cycle;
\draw (0,0) -- (2.4,1.6) -- (4.4,3.6) -- (6,6);
\draw (6,0) -- (2.4,1.6) -- (0,6) -- (4.4,3.6) -- cycle;
\draw[fill=black] (0,0) circle (1pt);
\draw[fill=black] (0,6) circle (1pt);
\draw[fill=black] (6,0) circle (1pt);
\draw[fill=black] (6,6) circle (1pt);
\draw[fill=black] (2.4,1.6) circle (1pt);
\draw[fill=black] (4.4,3.6) circle (1pt);
 \draw [fill=black] (4,2) circle (2pt);
 \draw [fill=black] (2.6,3.4) circle (2pt);
 \draw[dotted] (4,2) -- (2.6,3.4);
        \end{tikzpicture}   
        }
    \end{subfigure}
\\[10pt]
    \begin{subfigure}[t]{.25\linewidth}
    \subcaptionbox{Trianglicide. \label{fig:ABCDEF-T}}
    {
        \begin{tikzpicture}[scale=.4]
        \draw (0,0) -- (0,6) -- (6,6) -- (6,0) -- cycle;
\draw (0,0) -- (2.4,1.6) -- (4.4,3.6) -- (6,6);
\draw (6,0) -- (2.4,1.6) -- (0,6) -- (4.4,3.6) -- cycle;
\draw[fill=black] (0,0) circle (1pt);
\draw[fill=black] (0,6) circle (1pt);
\draw[fill=black] (6,0) circle (1pt);
\draw[fill=black] (6,6) circle (1pt);
\draw[fill=black] (2.4,1.6) circle (1pt);
\draw[fill=black] (4.4,3.6) circle (1pt);
 \draw[fill] (3,.5) circle (2pt);
 \draw[fill] (4.2,1.8) circle (2pt);
 \draw[fill] (5.4,2.8) circle (2pt);
 \draw[fill] (1,2.2) circle (2pt);
 \draw[fill] (2.6,3.4) circle (2pt);
 \draw[fill] (3.8,5) circle (2pt);
        \end{tikzpicture} 
        }
    \end{subfigure}
\qquad
    \begin{subfigure}[t]{.25\linewidth}
    \subcaptionbox{Our favorite: the $ACE$ example. \label{fig:ACE-T}}
    {
        \begin{tikzpicture}[scale=.4]
        \draw (0,0) -- (0,6) -- (6,6) -- (6,0) -- cycle;
\draw (0,0) -- (2.4,1.6) -- (4.4,3.6) -- (6,6);
\draw (6,0) -- (2.4,1.6) -- (0,6) -- (4.4,3.6) -- cycle;
\draw[fill=black] (0,0) circle (1pt);
\draw[fill=black] (0,6) circle (1pt);
\draw[fill=black] (6,0) circle (1pt);
\draw[fill=black] (6,6) circle (1pt);
\draw[fill=black] (2.4,1.6) circle (1pt);
\draw[fill=black] (4.4,3.6) circle (1pt);
\draw[fill=black] (1,2.2) circle (2pt);
\draw[fill=black] (4,2) circle (2pt);
\draw[fill=black] (3.8,5) circle (2pt);
%
        \end{tikzpicture}  
        }
    \end{subfigure}
\qquad
\caption{Some constrained triangulations}
\label{fig:examples-T}
\end{figure}
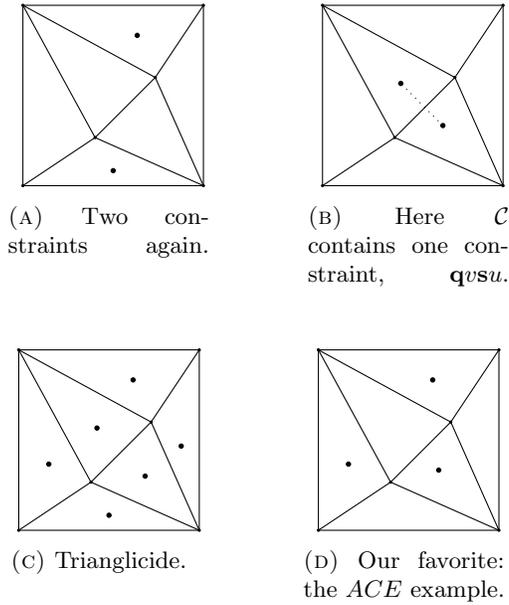

Figure \ref{fig:diag2} reproduces the honest triangulation $T$ of Example \ref{ex:2}.
Figure \ref{fig:examples-T} illustrates additional examples of the form $\T=(T,\CC)$, all with the same triangulation $T$.

\begin{example}[cf.~Example \ref{ex:2}]\label{ex:2AF}
If a constraint consists of the vertices of a single triangle, we indicate the constraint by marking the triangle.  For instance Figure \ref{fig:AF-T} has two constraints, each consisting of three vertices, indicated by the marks in the triangles.
\end{example}

\begin{example}[cf.~Examples \ref{ex:2}, \ref{ex:2BE}]
\label{ex:2BE(irr)}
Constraints consisting of vertices from multiple triangles are indicated by connecting the marks in the dual triangulation with dotted lines.  Figure \ref{fig:BE-T} has a single constraint $C$ consisting of the vertices of both marked triangles, i.e., the four vertices $\{\QQ,\SS,u,v\}$.  (Unlike with generalized dissections, these constraints do not come with a cyclic ordering because $T$ is given so we do not need to ensure that it is a disk.)  If we delete the edge $uv$, turning the two dead triangles into a quadrilateral, then we obtain the same 2-complex we get by poofing the dissection in Figure \ref{fig:BE-D}; the poofagon is the quadrilateral.  This figure shows one possibility for $\T(\Diss)$, where $\Diss$ is the (generalized) dissection of Figure \ref{fig:BE-D}.  The other is obtained by exchanging the edge $uv$ for the edge $\QQ\SS$.  

Incidentally, without the dotted line this example would be different; there would be two separate constraints that intersect in $u$ and $v$.  This is analyzed in Example \ref{ex:2BE(red)} in Section \ref{sec:combirr}.  This possibility is why we mark the triangles rather than shading them, as we did with poofagons.
\end{example}

\footnotetext{Ceci n'est pas une triangulation.}

\begin{example}[cf.~Example \ref{ex:2}]\label{ex:2(all)}
Figure \ref{fig:ABCDEF-T} shows an extreme example of a constrained triangulation.  Clearly this does not arise as $\T(\Diss)$ for any generalized dissection $\Diss$.
\end{example}

\begin{example}[$ACE$ again; cf.~Examples \ref{ex:2}, \ref{ex:2ACE-D}]
\label{ex:2ACE-T}
Figure \ref{fig:ACE-T} is the constrained triangulation for the $ACE$ example, so named because the living triangles are labeled $A,C,E$ in Figure \ref{fig:diag2}.  (Compare with Figure \ref{fig:ACE}.)  Here $\CC=\{\QQ vu,\RR \SS v,\SS\PP u\}$.  There is no way to realize this as a classical dissection, without killing one of the living triangles.  In Figure \ref{fig:ACE-D} the generalized dissection $\Diss$ has one upside-down triangle ($\SS uv$) and the same three constraints that are in $\CC$ (though there, technically, the constraints are cyclically ordered).  The constrained triangulation of Figure \ref{fig:ACE-T} is $\T(\Diss)$.
\end{example}

\begin{example}
If $\CC$ contains exactly one constraint and this constraint contains no more than 2 corners, then $\T$ is the type of object we studied in \cite{triangles1}.
Also there we usually required that the constraint be \emph{non-separating}.
\end{example}

\section{The space of drawings}\label{sec:X}

\subsection{Definitions and theorem}

We are now ready to introduce the space that allows us to talk about deforming a (generalized) dissection.

\begin{defn}\label{def:drawings} Let $\T=(T,\CC)$ be a constrained triangulation.  A \emph{drawing} of $\T$ is a map $\rho:\Vxs(T)\to\C^2$ such that
\begin{enumerate}
    \item for each $C\in\CC$ there is a line
    $\ell_C\subset\C^2$ such that 
    $\rho(v)\in\ell_C$ for each $v\in C$.
    \item the images of the corners form a parallelogram in $\C^2$; that is, $\rho(\PP)+\rho(\RR)=\rho(\QQ)+\rho(\SS)$;
\end{enumerate}
The space of all drawings, topologized as a subspace of $(\C^2)^{\Vxs(T)}$, is denoted $\dot X(\T)$.

A drawing is \emph{generic} if in addition to the above, we also have 
\begin{enumerate}
    \item[(3)] the $4$-gon (parallelogram) $(\rho(\PP),\rho(\QQ),\rho(\RR),\rho(\SS))$ is non-degenerate;
    \item[(4)] if $\{x,y,z\}$ is the vertex set of a living triangle in $\T$ then $(\rho(x),\rho(y),\rho(z))$ is a non-degenerate triangle in $\C^2$.
    \item[(5)] $\rho$ is injective (in particular this guarantees that the lines $\ell_C$ are uniquely defined); 
    \item[(6)] if $C,C'$ are distinct constraints with $C\cap C'\ne\emptyset$ then $\ell_C \ne \ell_{C'}$.  
\end{enumerate}
The closure in $\dot X(\T)$ (equivalently in $(\C^2)^{\Vxs(T)}$) of the set of all generic drawings of $\T$ is denoted $X(\T)$.  

We call $\T$ \emph{drawable} if there exists a generic drawing of $\T$, i.e., if $X(\T)$ is non-empty.
\end{defn}

The space $\dot X$ is evidently an algebraic variety, and as the closure of an open subset, $X$ consists of a union of components of $\dot X$.  In this section we give a parameterization of $X$ in the drawable case.  This shows that $X$ is rational and irreducible; it also follows that at most one component of $\dot X(\T)$ can contain a generic drawing.

\begin{thm}\label{thm:punt1}
Let $\T$ be drawable.  Then $X(\T)$ is an irreducible rational variety which is one of the components of $\dot X(\T)$.
\end{thm}

Some comments about the definition:

(1) Recall that we have defined the term ``generic'' already for generalized dissections.  
The connection is that if $\T$ is a (drawable) constrained triangulation, then every image of a generic drawing of $\T$ is a generic generalized dissection, and conversely, every generic generalized dissection $\Diss$ is the image of a generic drawing of the constrained triangulation $\T(\Diss)$. 

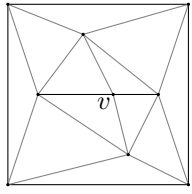
\begin{figure}
    \centering
    \begin{tikzpicture}[scale=.4]
        \draw (0,0) -- (0,6) -- (6,6) -- (6,0) -- cycle;
\draw [gray](1,3) -- (2.5,5) -- (5,3) -- (4,1) -- cycle;
\draw [gray](0,0) -- (1,3) -- (0,6) -- (2.5,5) -- (6,6) -- (5,3) -- (6,0) -- (4,1) -- cycle;
\draw (1,3) -- (5,3);
\draw [gray] (2.5,5) -- (3.5,3) -- (4,1);
\draw[fill=black] (0,0) circle (1pt);
\draw[fill=black] (0,6) circle (1pt);
\draw[fill=black] (6,0) circle (1pt);
\draw[fill=black] (6,6) circle (1pt);
\draw[fill=black] (1,3) circle (1pt);
\draw[fill=black] (2.5,5) circle (1pt);
\draw[fill=black] (5,3) circle (1pt);
\draw[fill=black] (4,1) circle (1pt);
\draw[fill=black] (3.5,3) circle (1pt);

\draw (3.2,2.7) node{$v$};
    \end{tikzpicture}
    \caption{This drawing is generic, even though three vertices are lined up horizontally.}
    \label{fig:generic}
\end{figure}

(2) Note also that this definition of generic is slightly more liberal than the usual concept of a \emph{general position} map of points into the plane (subject to (1) and (2) of course).  
Namely, we allow collections of vertices to be (accidentally) collinear, as long as such syzygies don't violate condition (6).
The dissection $D$ in Figure \ref{fig:generic} is generic, for example, and it is a generic drawing of $\T(D)$.  Compare with Figure \ref{fig:non-generic}.

(3) Observe that if $\Diss$ is a generic generalized dissection,  the slight ambiguity in defining $\T(\Diss)$ that arises in Lemma \ref{prop:DtoT} disappears in $X$, and so $$X(\Diss)=X(\T(\Diss))$$ is well-defined even though $\T(\Diss)$ isn't.  Likewise for $\dot X$.

(4) Examples of drawable triangulations include all honest triangulations ($\CC=\emptyset$) as well as any $\T(\Diss)$ for a generic generalized dissection $\Diss$.
(The latter follows from item (3) of Definition \ref{def:generalized}.)  

(5) Questions \ref{q:non-generic} and \ref{q:non-generic-classical} can be restated more precisely as follows.

\begin{questionnum}{1'}
Is $\T(\Diss)$ drawable for all (not necessarily generic) generalized dissections $\Diss$?
\end{questionnum}

\begin{questionnum}{2'}
Is $\T(D)$ drawable for all (not necessarily generic) dissections $D$?
\end{questionnum}

Our main interest here is (generalized) dissections, in which context Theorem \ref{thm:punt1} has the following consequence.

\begin{cor}\label{cor:punt1}
For any generalized dissection $\Diss$, the deformation space $X(\Diss)$ is either empty or a rational variety that is a single irreducible component of $\dot X(\Diss)$.
\end{cor}

\begin{proof}
If $X(\Diss)\ne\emptyset$ then $\T(\Diss)$ is drawable and Theorem \ref{thm:punt1} applies.
\end{proof}

\subsection{Combinatorial irreducibility and drawing orders}

We introduce some terminology before proving Theorem \ref{thm:punt1}.
Here is a pop quiz:  if $w,x,y,z$ are points in the plane, and $w,x,y$ are collinear, and $x,y,z$ are collinear, then must $w,x,y,z$ all be collinear?  The answer is no.  If $x\ne y$ then $x$ and $y$ determine a unique line and $w$ and $z$ must be on it.  But if $x=y$ then $w$ and $z$ can be anywhere.

\begin{figure}
    \centering
    \begin{tikzpicture}[scale=.5]
        \draw (0,0) -- (0,6) -- (6,6) -- (6,0) -- cycle;
\draw (0,0) -- (2.4,1.6) -- (4.4,3.6) -- (6,6);
\draw (6,0) -- (2.4,1.6) -- (0,6) -- (4.4,3.6) -- cycle;
\draw[fill=black] (0,0) circle (1pt);
\draw[fill=black] (0,6) circle (1pt);
\draw[fill=black] (6,0) circle (1pt);
\draw[fill=black] (6,6) circle (1pt);
\draw[fill=black] (2.4,1.6) circle (1pt);
\draw[fill=black] (4.4,3.6) circle (1pt);
 \draw [fill=black] (4,2) circle (2pt);
 \draw [fill=black] (2.6,3.4) circle (2pt);
    \end{tikzpicture}
    \caption{A combinatorially reducible $\T$.}
    \label{fig:BE-reducible}
\end{figure}
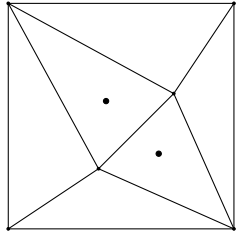

\begin{example}[cf.~Examples \ref{ex:2}, \ref{ex:2BE}, \ref{ex:2BE(irr)}]\label{ex:2BE(red)a}
Consider the constrained triangulation shown in Figure \ref{fig:BE-reducible}.  This has two separate constraints $C,C'$ (the marks are not joined by a dotted line).  Note that it has no generic drawings, because if we label the interior vertices $x$ and $y$, then by the pop quiz any drawing $\rho$ either has $\rho(x)=\rho(y)$ (violating condition (5) of genericity) or $\ell_C=\ell_{C'}$ (violating condition (6) of genericity).
\end{example}

This leads to the following definition and lemma, the proof of which is no different in the general case than it is in the above example.

\begin{defn} \label{def:combirr}
A constrained triangulation $\T=(T,\CC)$ is \emph{combinatorially irreducible} if there is at most one vertex in the intersection $C\cap C'$, for any two constraints $C,C'\in\CC$.
\end{defn}

\begin{lem}\label{lem:combirr}
Every drawable $\T$ is combinatorially irreducible.
\end{lem}

Our proof of Theorem \ref{thm:punt1} gives an explicit rational parameterization of $X(\T)$ in the drawable case.  
The tool we use is called a \emph{drawing order}.

Let $\T=(T,\CC)$ be drawable.  
It can nevertheless be difficult to actually draw $\T$, if one chooses an unfavorable order in which to draw the vertices. 
Let $\le$ be a total order on $\Vxs(T)$, with 
\begin{equation}\label{cornersfirst}
    \PP\le \QQ\le \SS\le \RR\le v \mbox{ for all interior }v\in\Vxs.
\end{equation}
Associated to $\le$ there is an integer-valued function $v\mapsto \a_v$ on $\Vxs$ defined as follows.
Label the vertices other than the corners by $v_1,\ldots,v_k$ so that $v_i\le v_j$ iff $i\le j$.  
Let $C\in\CC$ be a constraint.
For each $j=1,\dots,k$ define $C_{\le j}
=C\cap\{\PP,\QQ,\SS,\RR,v_1,\ldots,v_j\}\subset\Vxs(T)$, and say that $C$ is \emph{relevant} to $v_j$ if $v_j\in C$ and $|C_{\le j}| \ge 3$.
Now define $\a_\PP=\a_\QQ=\a_\SS=2$, $\a_\RR=0$, and for $j=1,\ldots,k$ 
\begin{equation}
    \a_j=\a_{v_j}:=2-\#\{C\in\CC \mid C \mbox{ is relevant to } v_j \}.
\end{equation}
We call the total order $\le$ a \emph{drawing order} if \eqref{cornersfirst} holds and also $\a_j\ge 0$ for each $j$.
Intuitively, we imagine trying to draw the vertices one by one in the order determined by $\le$.  When it is time to place $v_j$, the number of available degrees of freedom is (usually) $\a_j$.  As long as each $\a_j\ge0$, we will produce a drawing.

(It is not necessary to require that the corners come first, but it is convenient for the parameterization that follows.)

\begin{lem}\label{lem:drawingorder}
Every combinatorially irreducible $\T$ has a drawing order.
\end{lem}

\begin{proof}
Let $\T=(T,\CC)$ be combinatorially irreducible.
If $\Vxs(T)=\{\PP,\QQ,\RR,\SS\}$ then the order $\PP,\QQ,\SS,\RR$ is a drawing order.  So we may assume $T$ has interior vertices $v_1,\ldots,v_k$ with $k>0$. 

Let $R_1$ be the set of interior vertices of $T$ of valence less than $6$.  
This set is non-empty by an elementary argument about planar graphs, spelled out in Lemma \ref{lemma:redness}.  
For $j>1$, let $R_j$ be the set of interior vertices of the graph $G-\cup_{i<j}R_i$ that have valence less than $6$.  
(When a vertex is removed from a graph, any edges incident with the vertex are also removed.)  
Lemma \ref{lemma:redness} shows that each interior $v\in\Vxs(T)$ is in some $R_j$, because as long as $G-\cup_{i<j}R_i$ contains interior vertices of $T$, $R_j$ will be non-empty.

Now let $\le$ be any total order on $\Vxs$ such that \eqref{cornersfirst} holds and for all interior $v\in R_i, w\in R_j,$ if $i>j$ then $v<w$.
(The ordering within a given $R_j$ doesn't matter.)
We claim this is a drawing order.
The reason is that by construction each interior vertex $v$ is adjacent to fewer than six vertices that precede it in the order.  
By combinatorial irreducibility, any two constraints containing $v$ are otherwise disjoint.
So there must be fewer than three constraints that are relevant to $v$, because any such constraint would contribute at least two distinct neighbors of $v$ from among the earlier vertices.
Thus $\a_j\ge 0$ as desired.
\end{proof}

\begin{lem}\label{lemma:redness}
Let $G$ be a finite simple graph embedded in the plane such that the exterior face is bounded by the quadrilateral $\PP\QQ\RR\SS$.  
Label the other vertices $v_1, \ldots, v_k$ and assume $k>0$.
Then some $v_i$ has valence less than 6.
\end{lem}

\begin{proof}
Let $V=k+4$ denote the total number of vertices,
let $E$ denote the number of edges, and let $F$ denote the number of faces determined by the embedding of $G$ in the plane.
It is an easy consequence of Euler's equation $V-E+F=2$ that a finite simple planar graph must contain a vertex of valence less than 6.  
The point here is that $\PP,\QQ,\RR,\SS$ cannot be the only vertices with this property.
Our proof also involves nothing more than Euler's formula.

We may assume $G$ is connected, for otherwise we can find the vertex we seek in any connected component not touching the boundary.

Suppose for contradiction that each $v_i$ has valence at least 6.  
Then $G$ connected and $k>0$ imply that the sum of the valences of $\PP,\QQ,\RR,\SS$ is at least 9.
Summing valences we have $2E=\sum_v \valence(v) \ge 6k+9 =6(V-4)+9=6V-15$ and since $E\in\Z$ this implies $E\ge3V-7$.

Looking at faces we also have $2E=\sum_f \length(f) \ge 3F+1$ (where $\length(f)$ denotes the number of edges traversed by a path tracing out the boundary of the face $f$).

Now by Euler we have 
$$E+2=V+F\le \frac {E+7}{3} + \frac{2E- 1}{3}=E+2.$$
We must therefore have equality, so in particular all interior faces are triangles, which is the key to the rest of the argument.
This implies $\sum_{v\in\{\PP,\QQ,\RR,\SS\}} \valence(v) \ge 10$, whence $\sum_{v\in\{\PP,\QQ,\RR,\SS\}} \valence(v)=10$ (because otherwise the first inequality above would become $E\ge 3V-6$, a contradiction).
This in turn means only two (half-)edges emanate from the corners and go to the interior.
Say one emanates from $\PP$.  
Then the other cannot emanate from $\PP$ or $\QQ$ or $\SS$ since all interior faces are triangles.
Thus the other emanates from $\RR$.
Moreover these two edges must be the same edge $\PP\RR$, again because all interior faces are triangles.
Now no interior vertices can be present, meaning $k=0$, a contradiction.
\end{proof}

\subsection{Proof of theorem}

\begin{lem}\label{param}
Let $\T$ be drawable, and let $\le$ be a drawing order.  
There is a rational map 
$$g_{\leq}:\prod \C^{\a_v}\dashrightarrow X(\T)$$ 
parameterizing $X(\T)$.
Moreover, if $\rho$ is a generic drawing of $\T$ then $\rho$ has a unique preimage under $g_\le$, and in fact $g_\le$ is injective in a neighborhood this preimage. 
\end{lem}

Here we interpret $\C^0$ as the singleton $\{0\}$.

\begin{proof}
We define the rational map $$g:\prod\limits_{v\in\Vxs(T)} \C^{\a_v} \dashrightarrow (\C^2)^{\Vxs(T)}$$ as follows.
Let $x=(x_v)_{v\in\Vxs(T)}$ be coordinates on the domain of $g$.
The coordinate functions $g_v(x)$ of the point $g(x)$ are constructed inductively according to the chosen drawing order, as follows.

Recall that $\a_\PP=\a_\QQ=\a_\SS=2$ and $\a_\RR=0$.
We start by defining $g_\PP(x)=x_\PP\in\C^2$, $g_\QQ(x)=x_\QQ\in\C^2$, $g_\SS=x_\SS\in\C^2$, and $g_\RR(x)=x_\QQ+x_\SS-x_\PP$.
Then for each interior $v\in\Vxs(T)$, we assume that the coordinate functions $g_w(x)$ for vertices $w$ with $w<v$ have been defined.
To define $g_v$, we distinguish the three cases: $\a_v=0,1,2.$  

If $\a_v=2$, then $x_v$ is a point in $\C^2$, and we set $g_v(x)=x_v$.  

If $\a_v=1$, then $x_v\in\C$ is a number and there is a (unique) constraint that is relevant to $v$. We denote by   $y$ and $z$ the first two points with respect to $\le$ of this constraint, and we set $g_v(x)=x_v g_y(x) + (1-x_v) g_z(x)$. 

Finally if $\a_v=0$ then ($x_v=0$ and) there are two constraints $C$ and $C'$ relevant to $v$.
Let $y, z$ be the first two elements of $C$ (with respect to $\leq$), and let $y',z'$ be the first two elements of $C'$. 
Since $y,z,y',z' < v$, the rational functions $g_y(x),g_z(x),g_{y'}(x),g_{z'}(x)$ are already defined. 
Let $g_v(x)$ be the rational function in  $g_y(x),g_z(x),g_{y'}(x),g_{z'}(x)$ expressing the coordinates of the intersection of the line $L$ through $g_y(x),g_z(x)$ with the line $L'$ through $g_{y'}(x),g_{z'}(x)$.
This rational function can be computed explicitly, e.g., using Cramer's rule.  
There is a denominator.  
But $\T$ is assumed to be drawable, and in a generic drawing $\rho$ the two lines $L$ and $L'$ are distinct and non-parallel.
Thus this denominator is not identically zero, and so $g_v(x)$ is indeed a rational function.

We are now finished defining $g=g_\le$, and it remains to analyze its image.
Note that the domain of $g$ is irreducible, so the closure of $\im(g)$ is an irreducible algebraic variety.  We claim that this closure is exactly $X$.

Let $\rho$ be a generic drawing of $\T$.  
We define parameters $x=(x_v)$ (in order) such that $g(x)=\rho$.  If $\a_v=2$ the parameter is just $\rho(v)$. If $\a_v=1$ then by the injectivity of $\rho$ we know that $\rho(v)$ is an affine combination of the first two vertices in the relevant constraint for $v$, so the parameter $x_v$ is uniquely determined. If $\a_v=0$ then the fact that the lines $L$ and $L'$ are distinct and meet at the point $\rho(v)$ means that $\rho(v)$ is the (unique and correct) point parameterized by $g$.
Therefore, the image of $g$ contains all generic drawings, so the closure of $\im(g)$ contains $X$.

Moreover, if $g(x)=\rho$ is a generic drawing then there is an open neighborhood $U$ of $x=(x_v)\in\prod \C^{\a_v}$ such that $g(x')$ is a generic drawing for any $x'\in U$.
This is because conditions (1) and (2) of the definition of drawing are enforced by the definition of $g$, whereas (3), (4), (5), and (6) are open conditions.  
Thus $U$ maps into $X$.  Since $U$ is dense in the domain of $g$ and $X$ is closed, it follows that the closure of $\im(g)$ is contained in $X$, as desired.

Finally, we note that for distinct points $x',x''\in U$, if $v$ is the first vertex for which $x'_v\ne x''_v$ then $v$ has different images in $\C^2$ under the maps $g(x')$ and $g(x'')$.
So $g$ is injective on $U$.
\end{proof}

\begin{proof}[Proof of Theorem \ref{thm:punt1}]
The preceding lemma provides the necessary parameterization of $X(\T)$.  The inverse of $g$ is an algebraic map, so $X$ is indeed rational.  Moreover $X\subset \dot X$, $X$ is irreducible, and $X$ contains an open set of $\dot X$, so $X$ must be an irreducible component of $\dot X$.
\end{proof}

\begin{cor}\label{cor:Xdim}
Let $\T$ be drawable and let $\le$ be a drawing order.  Then $\sum_v \a_v$ is independent of choice of $\le$ and is equal to the dimension of $X(\T)$.
\end{cor}

Note that the dimension of $X$ agrees with the heuristic count, namely $6$ for the corners plus $2$ for each interior vertex minus $1$ for each vertex beyond the second in any constraint.

It is also worth noting that the affine group $\Aff=\Aff_2(\C)$ acts on $X$, and generic drawings have trivial stabilizers.
In particular $X(\T)$ is topologically a product of $\C^2$ (for translations) and a cone (for scaling) and is therefore contractible (if it is non-empty).
We do not know if the quotient $X/\Aff$ is contractible.

\subsection{Home field advantage}\label{sec:really}

    The space $\dot{X}$ consists of maps to $\C^2$. We conclude this section by arguing that from the point of view of drawing pictures, $\R^2$ would work just as well.  It may be interesting to study drawing spaces over other fields.

    \begin{defn}
    A constrained triangulation $\T$ is \emph{really drawable} if there is a generic drawing $\rho$ that maps all vertices into $\R^2$.
    Such a $\rho$ is called a \emph{real drawing}. 
    
    A constrained triangulation $\T$ is \emph{positively drawable} if there is a real drawing $\rho$ such that for any (oriented) triangle $(p,q,r)$ of $T$, the triangle $(\rho(p),\rho(q),\rho(r))$ in $\R^2$ is oriented positively.  
    Such a $\rho$ is called a \emph{positive drawing}.
    \end{defn}
    
    Here are some remarks about these definitions.
    \begin{enumerate}
        
        \item There certainly exist constrained triangulations that are positively drawable.  For instance every honest (unconstrained) triangulation is positively drawable.  Also, $\T(D)$ is positively drawable for any generic dissection $D$, since $D$ is a positive drawing of $\T(D)$. Conversely (the image of) any positive drawing of any $\T$ is a generic dissection.

        \item There exist (really) drawable constrained triangulations that are not positively drawable.  The smallest one is the $ACE$ example, shown in Figure \ref{fig:ACE-T}.
        
        \item If a constrained triangulation is drawable then it is really drawable. 
        To see this, find a drawing order and choose real parameters.  Almost all choices are in the domain of the parameterization $g$, because any denominator that vanishes at all real points would also vanish at all complex points.  Almost all real parameters in the domain of $g$ have image a real drawing.
        
        The inflection points of a complex cubic curve can be used to generate a linear system that is drawable but not really, by the Sylvester-Gallai theorem (see e.g.~\cite{thebook}).
        Existence of such things also follows from Mn\"ev universality \cite{mnev,mnev2}.  However they do not arise from planar triangulations.
        
        \item There exist constrained triangulations that are not drawable, by Lemma \ref{lem:combirr}.  
        Slightly less trivially, there are combinatorially irreducible $\T$'s that are not drawable.  For instance the constraints could force the boundary parallelogram to be degenerate. This may be the only obstruction to drawability in the combinatorially irreducible case.  See also Question 1.
    \end{enumerate}

\section{Musings about $\dot X$}
\label{sec:rambling}

This is an article about (generalized) dissections.  The notion of a constrained triangulation allows us to study deformations of generalized dissections, and as we have already observed, constrained triangulations of the form $\T(\Diss)$ have certain pleasant properties: they are always combinatorially irreducible, for instance, and if Question 1 has an affirmative answer then they are always drawable too.  The space $X(\Diss)$ is, as we have shown, an irreducible algebraic variety.

The collection of constrained triangulations includes many other interesting objects, though, that may be worthy of study on their own merits.  We conclude Part 1 by highlighting some examples and general questions about their drawing spaces, as well as some parallels with the theory of realizations spaces for oriented matroids.  Nothing in this section is central to the paper, although it is not entirely irrelevant either.  The reader who is anxious to get to the area relations can safely proceed to Part 2.

\subsection{The boundary}

Because of our interest in Monsky's theorem, we have so far only discussed constrained triangulations with four corners (see Definition \ref{def:constrained}), and we have required drawings to realize the boundary as a parallelogram.  Some of the issues we want to mention in this section are particular to that case, but many are not.  For the rest of this section we use the notation $\T_{\pentagon}$ to indicate a constrained triangulation with arbitrary boundary, whereas $\T$ continues to denote a constrained triangulation with four corners.

The spaces $\dot X$ and $X$ (Definition \ref{def:drawings}) can be defined for arbitrary $\T_{\pentagon}$ with the adjustments that condition (2) should be ignored and condition (3) should require the boundary, whatever it is, to be drawn as a non-degenerate polygon.  Combinatorial irreducibility (Definition \ref{def:combirr}) applies to $\T_{\pentagon}$ without modification.

While we are at it we give one more definition.
Given $\T$, we have defined both arbitrary drawings and generic drawings (Definition \ref{def:drawings}).  An intermediate type of drawing is one which satisfies conditions (1)--(4) of these definitions; we call these \emph{life-preserving}, as living triangles of $T$ are required to be drawn non-degenerately.  The closure of the life-preserving drawings is denoted $\hat X$.
Obviously $$X\subset\hat X\subset \dot X,$$
and like $X$, the space $\hat X$ is the closure of an open subset of $\dot X$, hence is a union of components of $\dot X$.

For $\T_{\pentagon}$, we make the same definition, modifying conditions (2) and (3) as we did earlier to define $\dot X$ and $X$.

    Our feeling is that $X(\T)$ captures the intuitive idea of deforming a dissection.
    However we acknowledge that this is to some extent a matter of taste; any of $X,\hat X,\dot X$ could reasonably be thought of as a deformation space for $\T$.

    A notational aside:  The $\dot{\phantom X}$ in $\dot X$ is meant to evoke the constant map, which is an element of $\dot X$, while the $\hat{\phantom X}$ in $\hat X$ resembles a triangle\footnote{$\overset{\triangle} X$ doesn't typeset very nicely.} to remind that (most) functions in $\hat X$ are faithful on the living triangles.  
    The notation $X$ has no decoration because it is used the most.
   
In the next subsection we make some conjectures about $\hat X$.

\subsection{Combinatorial reductions and $\hat X$}\label{sec:combirr}

If $\T_{\pentagon}=(T,\CC)$ is not combinatorially irreducible, then it can be decomposed into combinatorially irreducible factors.
If $u,v\in C\cap  C'$ for distinct $C,C'\in\CC$ (and distinct $u,v\in\Vxs(T)$) then the \emph{reduction} of $\T_{\pentagon}=(T,\CC)$ results in two \emph{combinatorial factors} (or just \emph{factors}) $\T_{\pentagon}'$ and $\T_{\pentagon}''$, where:
\begin{itemize}
    \item 
    $\T_{\pentagon}'=(T',\CC')$ where $T'=T$ and $\CC'=\CC$ except that $C, C'$ have been replaced by their union;
    \item 
    $\T_{\pentagon}''=(T'',\CC'')$ where $T''$ is the result of identifying vertices $u,v$ of $T$, and $\CC''$ is adjusted accordingly (removing any resulting constraints of size less than 3).
\end{itemize}

For various reasons, the factors resulting from these operations are not always constrained triangulations, and even if they are, they may not be combinatorially irreducible.
Moreover $\T_{\pentagon}$ may have multiple reductions.
Nevertheless, recursively continuing this procedure to its conclusion eventually leads to a ``factorization'' of $\T_{\pentagon}$ into a collection of combinatorially irreducible factors.

\begin{figure}
    \centering
    \begin{tikzpicture}[scale=.4]
        \draw (0,0) -- (0,6) -- (6,6) -- (6,0) -- cycle;
\draw (0,0) -- (2.4,1.6) -- (4.4,3.6) -- (6,6);
\draw (6,0) -- (2.4,1.6) -- (0,6) -- (4.4,3.6) -- cycle;
\draw[fill=black] (0,0) circle (1pt);
\draw[fill=black] (0,6) circle (1pt);
\draw[fill=black] (6,0) circle (1pt);
\draw[fill=black] (6,6) circle (1pt);
\draw[fill=black] (2.4,1.6) circle (1pt);
\draw[fill=black] (4.4,3.6) circle (1pt);
 \draw [fill=black] (4,2) circle (2pt);
 \draw [fill=black] (2.6,3.4) circle (2pt);
    \end{tikzpicture}
    \qquad
    \begin{tikzpicture}[scale=.3]
        \draw (0,0) -- (0,6) -- (6,6) -- (6,0) -- cycle;
\draw (0,0) -- (2.4,1.6) -- (4.4,3.6) -- (6,6);
\draw (6,0) -- (2.4,1.6) -- (0,6) -- (4.4,3.6) -- cycle;
\draw[fill=black] (0,0) circle (1pt);
\draw[fill=black] (0,6) circle (1pt);
\draw[fill=black] (6,0) circle (1pt);
\draw[fill=black] (6,6) circle (1pt);
\draw[fill=black] (2.4,1.6) circle (1pt);
\draw[fill=black] (4.4,3.6) circle (1pt);
 \draw [fill=black] (4,2) circle (2pt);
 \draw [fill=black] (2.6,3.4) circle (2pt);
 \draw[dotted] (4,2) -- (2.6,3.4);
    \end{tikzpicture}
        \begin{tikzpicture}[scale=.3]
            \draw (0,0) circle (0pt);
            \draw (0,2.8) -- (1,2.8);
            \draw (0,3.2) -- (1,3.2);
        \end{tikzpicture}    
    \begin{tikzpicture}[scale=.3]
        \draw (0,0) -- (0,6) -- (6,6) -- (6,0) -- cycle;
\draw (6,0) -- (0,6) ;
\draw (0,0) -- (2,4);
\draw (4,2) -- (6,6);
\draw [fill] (0,0) circle (1pt);
\draw [fill] (0,6) circle (1pt);
\draw [fill] (6,6) circle (1pt);
\draw [fill] (6,0) circle (1pt);
\draw [fill] (4,2) circle (1pt);
\draw [fill] (2,4) circle (1pt);
    \end{tikzpicture}
    \qquad
    \begin{tikzpicture}[scale=.3]
        \draw (0,0) -- (0,6) -- (6,6) -- (6,0) -- cycle;
\draw (4,2.5) -- (0,0);
\draw (4,2.5) -- (6,0);
\draw (4,2.5) -- (0,6);
\draw (4,2.5) -- (6,6);
\draw [fill] (0,0) circle (1pt);
\draw [fill] (0,6) circle (1pt);
\draw [fill] (6,6) circle (1pt);
\draw [fill] (6,0) circle (1pt);
\draw [fill] (4,2.5) circle (1pt);
    \end{tikzpicture}
    \caption{A combinatorially reducible $\T$ and its (two) irreducible factors, $\T'$ and $\T''$.}
    \label{fig:BE-reduction}
\end{figure}
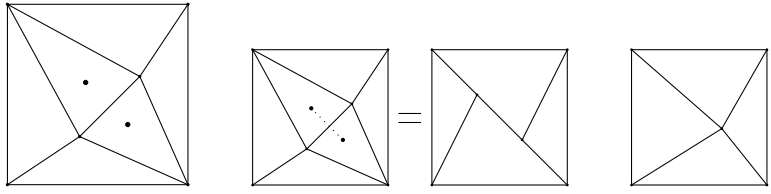

\begin{example}[cf.~Examples \ref{ex:2}, \ref{ex:2BE}, \ref{ex:2BE(irr)}, \ref{ex:2BE(red)a}]\label{ex:2BE(red)}
The basic example to keep in mind is the $\T$ shown in Figure \ref{fig:BE-reduction}.  All these figures have appeared before; this is a parallelogram.
There are two constraints in $\T$.  In Figure \ref{fig:BE-reduction} we have equated the factor $\T'$ with a generic drawing of it.  We have also shown $\T''$, which is an honest triangulation.  Both $\T'$ and $\T''$ are combinatorially irreducible and drawable; the figure shows generic drawings of both. Note that neither of these is considered a generic drawing of $\T$, though for different reasons.

In this case $\dot X(\T')=X(\T')$ (though this is not totally obvious) and $\dot X(\T'')=X(\T'')$.  Both are irreducible components of $\dot X(\T)$, which has no other components.  
\end{example}

This is a good time to point out that combinatorial irreducibility is not necessary for the existence of a drawing order (compare Lemma \ref{lem:drawingorder}).
Recall that from a drawing order $\le$ we produce a parameterization $g_{\le}$ of a component of $\dot X(\T)$.
In the combinatorially irreducible case any drawing order will yield the same component, namely $X(\T)$.
On the other hand in the current example, with the interior vertices labeled $u$ and $v$, the drawing order $\PP \QQ \SS \RR uv$ yields a parameterization of the component $X(\T'')$, whereas the drawing order\footnote{This is technically not a drawing order but the point remains.} $\PP \QQ uv \SS \RR$ leads to a parameterization of the other component $X(\T')$.

\begin{conjecture}\label{conj:irr}
If $\T_{\pentagon}$ is a constrained triangulation with arbitrary boundary polygon, then $\T_{\pentagon}$ is combinatorially irreducible if and only if $\hat X(\T_{\pentagon})$ is an irreducible variety.
\end{conjecture}

In the case we have focused on for the majority of this paper, i.e., constrained triangulations $\T$ with the boundary condition, an extra condition is required to make the analogous conjecture possible.

    \begin{defn} 
    A constrained triangulation $\T$ (of a parallelogram) is \emph{toroidally irreducible} if (a) it is combinatorially irreducible and (b) there do not exist constraints $C,C'$ with $\{\PP,\QQ\} \subset C$ and $\{\RR,\SS\} \subset C'$ and (c) there do not exist constraints $C,C'$ with $\{\QQ,\RR\} \subset C$ and $\{\PP,\SS\} \subset C'$.
    \end{defn}
    
    This is sort of like saying that $\T$ is combinatorially irreducible after identifying opposite edges of the boundary to make $T$ into (a triangulation of) a torus.
    We will not spell out the reduction process but one can imagine that $\T'$ has a (single) constraint that ``wraps around'' the torus, and $\T''$ is a ``constrained triangulation of a segment.''
    
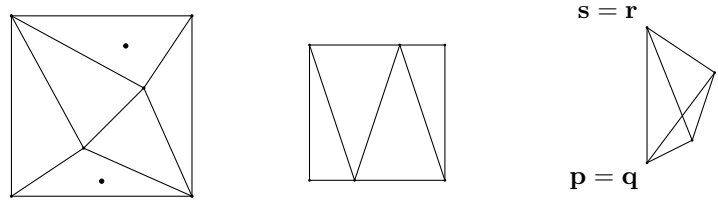
\begin{figure}
    \centering
    \begin{tikzpicture}[scale=.4]
        \draw (0,0) -- (0,6) -- (6,6) -- (6,0) -- cycle;
\draw (0,0) -- (2.4,1.6) -- (4.4,3.6) -- (6,6);
\draw (6,0) -- (2.4,1.6) -- (0,6) -- (4.4,3.6) -- cycle;
\draw[fill=black] (0,0) circle (1pt);
\draw[fill=black] (0,6) circle (1pt);
\draw[fill=black] (6,0) circle (1pt);
\draw[fill=black] (6,6) circle (1pt);
\draw[fill=black] (2.4,1.6) circle (1pt);
\draw[fill=black] (4.4,3.6) circle (1pt);
 \draw[fill] (3,.5) circle (2pt);
 \draw[fill] (3.8,5) circle (2pt);
    \end{tikzpicture}
    \qquad \qquad
    \begin{tikzpicture}[scale=.3]
        \draw (6,6) -- (6,0) -- (0,0) -- (0,6) -- cycle;
\draw (0,6) -- (2,0) -- (4,6) -- (6,0);
\draw[fill=black] (6,6) circle (1pt);
\draw[fill=black] (6,0) circle (1pt);
\draw[fill=black] (0,0) circle (1pt);
\draw[fill=black] (0,6) circle (1pt);
\draw[fill=black] (2,0) circle (1pt);
\draw[fill=black] (4,6) circle (1pt);
\draw (0,-.75) circle (0pt);
    \end{tikzpicture}
    \qquad \qquad
    \begin{tikzpicture}[scale=.3]
        \draw (0,0) -- (0,6);
\draw (0,0) -- (2,1) -- (3,4) -- (0,6);
\draw (2,1) -- (0,6);
\draw (0,0) -- (3,4);
\draw[fill=black] (0,0) circle (1pt);
\draw[fill=black] (0,6) circle (1pt);
\draw[fill=black] (2,1) circle (1pt);
\draw[fill=black] (3,4) circle (1pt);
\draw (0,0) node[anchor=north east]{$\PP=\QQ$};
\draw (0,6) node[anchor=south east]{$\SS=\RR$};
    \end{tikzpicture}
    \caption{A toroidally reducible $\T$ (left) and drawings of its two irreducible factors, $\T'$ and $\T''$.}
    \label{fig:AF-reduction}
\end{figure}

    Figure \ref{fig:AF-reduction} (left) exhibits toroidal reducibility.  Here $\CC$ has two constraints, $\PP \QQ u$ and $\RR \SS v$ (using our usual notation).  This is combinatorially irreducible but not toroidally irreducible.  
    
    The space $\dot X(\T)$ has two components.  One component is $X(\T)$, consisting of drawings $\rho$ with non-degenerate boundary $\PP \QQ\RR\SS$ and with $u$ on the line $\PP \QQ$ and $v$ on the line $\RR\SS$.
    One such drawing is shown in Figure \ref{fig:AF-reduction} (middle); these drawings are (almost all) generic.  This component coincides with $X(\T')$ and $\dot X(\T')$.  
    The other component of $\dot X$ consists entirely of non-generic drawings $\rho$ having $\rho(\PP)=\rho(\QQ)$ and $\rho(\RR)=\rho(\SS)$; these might be thought of as ``dissections of a segment.''
    Both components are 8-dimensional (2-dimensional after quotienting by the affine group action).
    
\begin{conjecture}\label{conj:toroidal}
If $\T$ is toroidally irreducible then $\hat X(\T)$ is irreducible. 
\end{conjecture}

\subsection{Components of $\dot X$}

    We give some more examples to illustrate the differences between $X,\hat X,$ and $\dot X$.
    
    For honest triangulations we have $X=\hat X=\dot X$; this space is non-empty and irreducible and isomorphic to an affine space. 
    
    Example \ref{ex:2BE(red)} (Figure \ref{fig:BE-reduction}) has $\emptyset=X \ne \hat X=\dot X$, and the latter has two components of the same dimension.  This example is combinatorially reducible, hence not (generically) drawable.  The components of $\hat X=\dot X$ are $X(\T')$ and $X(\T'')$ where $\T'$ and $\T''$ are the  factors of $\T$.  In other words although there is no generic drawing of $\T$, every drawing of $\T$ is close to a generic drawing of one of the two combinatorially irreducible factors of $\T$.
    
    Example \ref{ex:2AF} (Figure \ref{fig:AF-T}) is drawable, and we have $\emptyset\ne X\ne\hat X = \dot X$; again $\dot X$ has two components of the same dimension.  This example is not toroidally irreducible.  The components of $\hat X=\dot X$ are again $X(\T')$ and $X(\T'')$ where the two factors are obtained by the toroidal reduction alluded to above. The first of these is also $X(\T)$.  

\begin{figure}
        \centering
        \begin{tikzpicture}[scale=.4]
            \draw (0,0) -- (0,6) -- (6,6) -- (6,0) -- cycle;
\draw (4,2.5) -- (0,0);
\draw (4,2.5) -- (6,0);
\draw (4,2.5) -- (0,6);
\draw (4,2.5) -- (6,6);
\draw [fill] (0,0) circle (1pt);
\draw [fill] (0,6) circle (1pt);
\draw [fill] (6,6) circle (1pt);
\draw [fill] (6,0) circle (1pt);
\draw [fill] (4,2.5) circle (1pt);
\draw [fill] (3.5,1) circle (2pt);
\draw [fill] (5,2.8) circle (2pt);
\draw [fill] (1.5,2.5) circle (2pt);
        \end{tikzpicture}
        \caption{Collateral damage.}
        \label{fig:collateral}
    \end{figure}
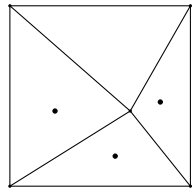

    A new example shown in Figure \ref{fig:collateral} exhibits $\emptyset=X=\hat X\ne \dot X$, and $\dot X$ has just one component, consisting of drawings with all five vertices on a line. The phenomenon on display here is called \emph{collateral damage}.
    
    At the opposite extreme from the honest case, suppose that $\T=(T,\CC)$ where the vertices of each triangle of $T$ form a constraint. Here of course $\hat X(\T)=X(\T)=\emptyset$.  The space $\dot X$ has a component consisting of drawings in which all points are collinear, but there may also be other components of smaller dimension.  (In Example \ref{ex:2(all)}, $\dot X$ has just one component.)  The space $\dot X(\T)$ plays a role in our study because it is a model for the base locus of the area map $\Area:X(T)\dashrightarrow Y(T)$ associated to the honest $T$.  (In fact that base locus is always contained in $\dot X(\T)$, though this may be a proper containment.)  We analyze this base locus in a forthcoming paper.
    
    We do not know if it is possible to have $X\ne\hat X\ne \dot X$.  What about in the drawable case, i.e., is $\emptyset \ne X\ne\hat X\ne \dot X$ possible?

    It may be the case that components of $\dot X(\T)$ can always be interpreted as $X(\T^*)$ for the various combinatorial factors $\T^*$ of $\T$.  Some of the components, those making up $\hat X$, are $X(\T^*)$ for the factors that are themselves drawable constrained triangulations.  If there is only one of these with non-degenerate boundary, it is also $X(\T)$.
    However, this picture is merely conjectural.

    \begin{question}\label{q:dimX}
    How many components does $\dot X$ have, and what are their dimensions? 
    \end{question}
    
    In particular we do not know if $\dot X$ can have components of dimension larger than $\dim X$ or $\dim \hat X$, when either of these two is non-empty.  
    (In the drawable case, of course, Corollary \ref{cor:Xdim} gives the dimension of $X$.)
    As this may involve subtle issues in incidence geometry, the answer may again (like Question 1) depend on the underlying field.  

    Recall that the heuristic dimension count is $6$ for the corners, $2$ for each additional vertex, and $-1$ for each vertex beyond the second in any individual constraint.  (For $\T_{\pentagon}$, the heuristic is the same except the boundary contributes $2m$ if it has $m$ vertices.)  Corollary \ref{cor:Xdim} verifies this for the component $X$ of $\dot X$, when the former is non-empty.  
    It is possible that in the absence of collateral damage this holds for all $\T$, and even $\T_{\pentagon}$.
    As far as we know, though, the dimension could be higher than the heuristic indicates, because the constraints could be redundant (in obvious or subtle ways).  
    Here is an ``obvious'' example:  
    constraints $abc,bcd,acd$ are equivalent to $abcd$.  
    Thus these three really only cut the dimension down by 2.  
    This is because the third constraint selects a component of the reducible variety determined by the first two.  
    The heuristic is too low by 1.
    
    A non-obvious redundancy could arise if there were a (non-trivial) incidence theorem, such as Pappus.  
    Here there are nine points, and the collinearity of eight specific triples implies the collinearity of a ninth.  
    So, after accounting for the eight hypothesized constraints, further including the ninth lowers the heuristic dimension count but doesn't actually change the variety.  
    
    Things like this (probably) are what make the dimension of the realization space algorithmically intractable for general point/line configurations.  (See below.)
    Fortunately, Pappus' theorem does not come into play for us, because we only work with triangulations of a disk whereas the configuration of Pappus' theorem is non-planar.\footnote{If each triple that is a collinearity hypothesis of Pappus' theorem is made into a triangle, then the resulting 2-complex made of eight triangles does not embed in the plane, because its 1-skeleton contains (a subdivision of) the graph $K_{3,3}$.} 
    However there may be other incidence theorems and we do not know whether any non-obvious redundancies can arise in the setting of constrained triangulations.

\subsection{Realization spaces of oriented matroids}

    The issues we have mentioned so far in this section are reminiscent of general questions about realizing configurations of points and lines.
    We now focus on this analogy, and we present a few variations on several problems about oriented matroids that are known to be difficult.

    At points of $\dot X$, although certain triangles are required to be degenerate, the others have no restriction one way or the other.  As a result $\dot X$ contains all constant functions, and $\dot{X}$ decomposes as a product of $\C^2$ (due to translations) and a cone (due to scaling).  In particular $\dot{X}$ is contractible.  It may be interesting to study the topology of the quotient of $\dot X$ by the affine action.
    
    By contrast, in the context of point/line configurations (commonly described in the language of oriented matroids), there are point/line configurations with disconnected realization space.  
    This is the ``isotopy problem'' for point/line configurations, solved in the 1980's by various people including the Mn\"ev universality theorem \cite{mnev,mnev2}.
    This means there are different drawings of points, with the same pre-specified incidence relations, that are not isotopic to each other through such configurations. 
    This is a restricted sense of isotopy, though, as in a realization space one is not allowed to introduce degeneracies (even temporarily) along the way.  In this sense our deformation spaces are fundamentally different.

    If Conjecture \ref{conj:irr} is true, it would suggest that these ``planar'' systems are significantly simpler than general systems, just as planar graphs are significantly simpler than general graphs.

    Nevertheless we suspect that Question \ref{q:dimX} and its relatives are difficult even for constrained triangulations.
        
    It is worth writing down the following (probably intractable but more basic) question.
    
    \begin{question}
    Given a finite set $\Vxs$ and a finite collection of subsets $\CC$ of $V$, what is the dimension of 
    the space of those maps $\rho:V\to\C^2$ such that for each set $C\in\CC$, the set of points 
    $\{\rho(v) : v\in C\}$ lies in a line?
    \end{question}
    The same question can be asked with $\C^2$ replaced by $\F^n$ for any field $\F$.  
    We would be interested to know if there are results about the hardness of this problem.  
    
    A related question that we know virtually nothing about is the following.  Fix a finite simplicial complex $T$ and a number $n$ such that $T$ embeds in $\R^n$.  Let $d$ be a function on the simplices of $T$ such that $d(\sigma)\le \dim(\sigma)$ for all $\sigma$, and also $d(\sigma)\le d(\tau)$ if $\sigma\subset\tau$.  What is the nature of the space $X$ of maps $\rho:T\to\R^n$ satisfying $\dim(\rho(\sigma))=d(\sigma)$ for all simplices $\sigma$?  What about maps satisfying $\dim(\rho(\sigma))\le d(\sigma)$?

\part{Area relations}

\bigskip

We now shift gears and begin our study of Monsky's theorem and the polynomials $f$ and $p$ discussed in the introduction.
Our principal contribution is to extend Monsky's theorem to the deformation spaces $X$ that we defined in Part 1.
We then explore the consequences of this extension for $f$ and $p$.

Henceforth all constrained triangulations $\T$ will be assumed to have square boundary.

\section{Area of a triangle}

Let $\F$ be a field of characteristic not equal to 2, and let $p_i=(x_i,y_i)$ for $i=1,2,3$ be three points in the affine plane $\F^2$.  We define the area of the (ordered) 
triangle $\Delta=(p_1,p_2,p_3)$ to be 
$$\Area(\Delta)=\Area(p_1p_2p_3)=\frac 12
\left|
\begin{matrix}
1 & 1 & 1 \\
x_1 & x_2 & x_3 \\
y_1 & y_2 & y_3 
\end{matrix}
\right|.
$$
Note that if $\F=\R$ then this is the usual signed area function.  We will also use this definition when $\F$ is $\C$ or a function field.

Note also that regardless of the field, $\Area(p_1p_2p_3)=0$ if and only if $p_1,p_2,p_3$ lie on a line in $\F^2$.

\section{Monsky Theorems}

If $D$ is a (classical) dissection of the unit square into triangles with areas $a_1,\ldots,a_n$, then Monsky's theorem gives a polynomial $f$ with integer coefficients such that 
\begin{equation}\label{eq:f-nonhomog}2f(a_1,\ldots,a_n)=1. \end{equation}  

In this section ,we use a modification of Monsky's argument, carried out over the field of rational functions in the vertex coordinates, to show that $f$ can be chosen to depend only on $\T(D)$ and not on $D$, meaning that the $a_i$ can represent the areas of the triangles in any drawing of $\T(D)$, including drawings that are not positive.
Accordingly, the equation \eqref{eq:f-nonhomog} needs to be modified to take into account the total area of the drawing; see \eqref{eq:f} below.
Moreover, because we carry out the argument in the setting of abstract triangulations, the theorem will apply equally well to generalized dissections as to dissections, even when the former have no positive drawings.  Equation \eqref{eq:f} will also hold for limits of drawings.

We give two versions of this argument, the first for honest triangulations in Section \ref{sec:+} and the second that incorporates the constraints in Section \ref{sec:++}.
In the presence of constraints, these results have significant computational benefit over the approach taken in \cite{triangles1}.  We discuss this further in Section \ref{sec:compute}.

\subsection{Monsky homogenized and deformed}

We have set up our drawing spaces so that the boundary can be mapped to an arbitrary parallelogram, rather than just the unit square.  We state our generalization of Monsky's Theorem in a similar spirit.  For this purpose, it makes sense to homogenize the equation of Monsky's Theorem.  This is quite easy to do.  Take a dissection of the unit square, and take any $f\in \Z[A_1, \dots, A_n]$ satisfying Monsky's theorem for this dissection.  Note that the polynomial $\sigma=A_1+\cdots+A_n$ evaluates to $1$ when the areas $a_i$ are plugged in. Thus if we homogenize $f$ with respect to the homogenizing variable $\sigma$ we obtain a homogeneous polynomial $\hat{f}$ satisfying
\begin{equation}\label{eq:f}
2\hat{f}(a_1, \dots a_n) = \sigma(a_1, \dots, a_n)^e,
\end{equation}
where $e$ is the degree of $f$. This  relation, now homogeneous, has the advantage of being affine invariant; that is, this equation will hold not only for the original dissection, but also for any affine image of it. 

We also wish to find relations that are invariant under deformations.  The following example illustrates this and introduces some notation used in the statement of the generalized Monsky Theorem.
\begin{example}[cf.~Examples \ref{ex:1}, \ref{ex:1-D}]\label{ex:diag1}
In Example \ref{ex:1}, we introduced the dissection $D$ of Figure \ref{fig:diag1-labeled}. 
\begin{figure}
\begin{center}
    \begin{tikzpicture}[scale=.5]
        \draw (0,0) -- (0,6) -- (6,6) -- (6,0) -- cycle;
\draw (4,2.5) -- (0,0);
\draw (4,2.5) -- (6,0);
\draw (4,2.5) -- (0,6);
\draw (4,2.5) -- (6,6);
\draw [fill] (0,0) circle (1pt);
\draw [fill] (0,6) circle (1pt);
\draw [fill] (6,6) circle (1pt);
\draw [fill] (6,0) circle (1pt);
\draw [fill] (4,2.5) circle (1pt);
\draw (3.5,1) node{$A$};
\draw (5,2.8) node{$B$};
\draw (3.5,4.5) node{$C$};
\draw (1.5,2.5) node{$D$};
    \end{tikzpicture}
    \caption{}
    \label{fig:diag1-labeled}
\end{center}
\end{figure}
If this is the unit square, and the central vertex has coordinates $(x,y)$, then the areas are
$\tilde{Z}_A = \frac12 y, \tilde{Z}_B = \frac12 (1-x), \tilde{Z}_C= \frac12 (1-y), \tilde{Z}_D = \frac12 x$.  The tildes over the $Z$'s indicate that the corners have been fixed to those of the unit square.  In a moment we will switch from $\tilde{Z}$ to $Z$, when we allow the boundary to be an arbitrary parallelogram.  In the spirit of finding relations among the areas that are preserved under deformations, we consider $\tilde{Z}_A, \tilde{Z}_B, \tilde{Z}_C, \tilde{Z}_D$ to be polynomials (or rational functions) living in the field of rational functions in the two coordinate variables $x$ and $y$. Contrast this with the $a_i$ of Monsky's Theorem, which are real numbers.  We seek algebraic relations among the four rational functions $\tilde{Z}_A, \tilde{Z}_B, \tilde{Z}_C, \tilde{Z}_D$.  

In this case, finding such relations is not hard.  Corresponding to the geometric observation that the bottom and top triangles add up to half the area of the square, algebraically we have
$2(\tilde{Z}_A+\tilde{Z}_C)=1$, or more homogeneously,
$$
2(\tilde{Z}_A+\tilde{Z}_C)= \sigma, 
$$
where $\sigma = \tilde{Z}_A+\tilde{Z}_B+\tilde{Z}_C+\tilde{Z}_D. $
Thus the polynomial $f(A,B,C,D)=A+C$ satisfies \eqref{eq:f-nonhomog}, provided the boundary is a unit square. 

Having found a homogeneous relation of the form $2f=\sigma^e$ among the $\tilde{Z}$'s, we now observe that the same relation holds even if the boundary is an arbitrary parallelogram.  Precisely, consider the rational function field $S$ in the eight variables $x_v, y_v$ where $v$ is one of the four vertices other than $\RR$, and inside $S$, define $x_{\RR}=x_\QQ+x_\SS-x_\PP$ and $y_\RR=y_\QQ+y_\SS-y_\PP$. Let $Z_A, Z_B, Z_C, Z_D\in S$ be the homogeneous quadratic polynomials in these eight variables expressing the areas of the triangles without fixing the corners.  Then by an easy argument invoking affine invariance, the identical relation $2f=\sigma^e$ will hold among the $Z$'s.  This argument is spelled out in detail in Corollary \ref{cor:fixcorn}.
\end{example}

\subsection{Honest triangulations}\label{sec:+}

We now give our first modification of Monsky's argument, designed for honest triangulations.

The clever uses of ultranorms and Sperner's lemma in the proofs trace directly to Monsky's original theorem \cite{monsky}.  The use of Sperner's lemma built on an earlier approach due to Thomas \cite{thomas}.
We emphasize that we are adapting those ideas to our current context.
We mimic the treatment in Pete Clark's class notes \cite{peteclark} which fleshes out some of the steps.

We first establish some notation. 
Let $T$ be a triangulation of a square with $k$ interior vertices and $n=2k+2$ triangles.
Let $S$ be the rational function field $\Q(x_v,y_v)$ for $v\in\Vxs(T)\setminus\{\RR\}$.  In $S$, set $x_{\RR}=x_\QQ+x_\SS-x_\PP$ and $y_\RR=y_\QQ+y_\SS-y_\PP$. The field $S$ is the coordinate function field of the space of drawings of $T$.

\begin{thm}[Monsky+] Let $T$ be a triangulation of a square and let $S$ be the coordinate function field of its drawing space. For each triangle $\Delta_j$, let $Z_j\in S$ be the homogeneous quadratic polynomial in $\{x_v,y_v\}$ expressing the area of $\Delta_j$.  Define $\sigma \in S$ by 
 $$\sum Z_j=\sigma.$$
Then there exists a homogeneous polynomial $f_T$ with integer coefficients in $n$ variables such that 
$$2f_T(Z_1,\ldots,Z_n)=\sigma^e$$
in $S$, for some non-negative integer $e$.
\end{thm}

\begin{defn}
If $T$ is a triangulation of a square with $n$ triangles, then any polynomial $f$ satisfying the conclusion of Theorem Monsky+ is called a \emph{Monsky polynomial} for $T$.
\end{defn}

For example, the polynomials $A+C, B+D,$ and $2A-B+2C-D$ are all Monsky polynomials for the example of Figure \ref{fig:diag1-labeled}.

\begin{proof}
Let $R$ be the subring $\Z[Z_1, \dots, Z_n]\subset S$.
Note that $\sigma\in R$, and we endeavor to show that $\sigma$ is an element of the ideal $\sqrt{(2)}$ of $R$.  It then follows that some power $\sigma^e$ may be written as 2 times an integer polynomial in the $Z_i$.  Furthermore, this polynomial can be chosen to be homogeneous of degree $e$ since both $\sigma$ and all the $Z_i$ are homogeneous of the same degree. 

Thus all is reduced to showing $\sigma\in\sqrt{(2)}$. Assume this is not the case.  Then there is a minimal prime $\pp$ containing $(2)$ such that $\sigma\not\in\pp$ (since $\sqrt{(2)}$ is the intersection of all prime ideals containing $(2)$).
By Krull's principal ideal theorem $\pp$ has height $1$.

Let $R_\pp$ be the ring $R$ localized at the prime ideal $\pp$ and $\bar R$ be the integral closure of $R_\pp$.   
The ideal $\pp'=\pp R_\pp$ in $R_\pp$ has height 1.  
Let $\qq$ be a prime ideal of $\bar R$ lying over $\pp'$.  
Then $\qq$ has height 1 in $\bar R$ \cite{irena}.

By the Mori-Nagata theorem $\bar R$ is a Krull domain, 
hence $\bar R_\qq$ is a discrete valuation ring containing $\bar R$.
The valuation on $\bar R_\qq$ yields a non-Archimedean ultranorm $\|\cdot\|$ on the fraction field of $\bar R_\qq$ (which is also the fraction field of $R$).  
Since $Z_j\in R\subset \bar R_\qq$, we have $\|Z_j\|\le 1$ for each $j$.  
In addition $2\in\pp$ so $2\in\qq\bar R_\qq$ so $\|2\|<1$.
Furthermore $\sigma\not\in\pp$ so $\sigma$ is a unit in $R_\pp$, hence in $\bar R$, hence $\|\sigma\|=1$.

Extend the ultranorm $\|\cdot\|$ to $S$, referring to  \cite{jahnke}
as necessary.  

Now, following Monsky, we color each point $\phi=(\phi_x,\phi_y)$ of the plane $S^2$ with one of the colors $A,B,C$ via the following comparisons:  

\begin{itemize}
    \item if $\|\phi_x\| \ge \|\phi_y\|$ and $\|\phi_x\| \ge \|1\|$ then $\phi$ gets color $A$;
    \item else if $\|\phi_y\| \ge \|1\|$ then $\phi$ gets color $B$;
    \item else $\phi$ gets color $C$.
\end{itemize}
In other words we color $A,B,C$ according to which of $\phi_x, \phi_y,$ or $1$ has the largest norm, breaking ties in that order.  (See \cite{monsky}.)

Monsky proved two lemmas which hold in this context exactly as he proved them, using the defining properties of the ultranorm, namely $\|\alpha\beta\|=\|\alpha\|\,\|\beta\|$ and $\|\alpha+\beta\|\le\max\{\|\alpha\|,\|\beta\|\}$, with equality if $\|\alpha\|\ne\|\beta\|$.
We leave the verifications as exercises.

\begin{lem}[Monsky]  
The color of $\phi$ agrees with the color of $\phi+\psi$ for any $C$-colored $\psi$.
\end{lem}
\begin{lem}[Monsky]  
Any triangle $\Delta$ whose vertices are colored $ABC$ satisfies $\|\Area(\Delta)\|>1$.
\end{lem}

We intend to use this coloring of $S\times S$ to induce a coloring of the vertices of $T$.  We do this as follows.

Let $M:S\times S\to S\times S$ be the unique affine transformation on $S\times S$ taking $(x_\PP,y_\PP)$ to $(0,0)$, $(x_\QQ,y_\QQ)$ to $(1,0)$, and $(x_\SS,y_\SS)$ to $(0,1)$.
Note that the determinant of $M$ is
$\left|
\begin{matrix}
x_\QQ-x_\PP & x_\SS-x_\PP \\
y_\QQ-y_\PP & y_\SS-y_\PP
\end{matrix}
\right|^{-1}$,
which equals the nonzero element $\frac 1{\sigma}$ of $S$ (so in particular $M$ exists).
Thus for any triangle $\Delta$ in $S\times S$ we have
$$\Area(M\Delta)=\frac 1{\sigma} \Area(\Delta)$$
and since $\|\sigma\|=1$,
$$ \|\Area(M\Delta)\|=\|\Area(\Delta)\|.$$
(Note that we are computing $\Area$ in the field $S$.)

Now, as promised, the coloring of $S\times S$ induces a coloring of the vertices of $T$ by assigning to the vertex $v$ of $T$ the color of the point $M(v_x,v_y)$.
This colors the corners $\PP \QQ\RR\SS$ with the colors $CAAB$ and is therefore a Sperner coloring, so there is an $ABC$ triangle $\Delta_j$.  
By Monsky's lemma 2, $\|\Area(M(\Delta_j))\|>1$.  
But $\|\Area(M(\Delta_j))\|=\|\Area(\Delta_j)\|=\|Z_j\|\le 1$, a contradiction.
\end{proof}

When working with the polynomial $f$, it is useful to exploit the fact that any parallelogram is affinely equivalent to the unit square, thereby allowing us to remove the $x,y$ variables corresponding to corners.  We make now this idea precise.

Given a triangulation $T$, consider the fixed-corner function field $\tilde{S}$, defined to be the rational function field in the $2k$ variables ${\tilde{X}_v, \tilde{Y}_v}$, where $v$ denotes an interior vertex. For corner vertices $v$, define elements ${\tilde{X}_v, \tilde{Y}_v}\in  \tilde{S}$, by 
$(\tilde{X}_{\PP}, \tilde{Y}_{\PP}) =  (0,0)$, $(\tilde{X}_{\QQ}, \tilde{Y}_{\QQ}) =  (1,0)$, $(\tilde{X}_{\RR}, \tilde{Y}_{\RR}) =  (1,1)$, $(\tilde{X}_{\SS},\tilde{Y}_{\SS}) =  (0,1)$.

\begin{cor}\label{cor:fixcorn} Let $T$ be a triangulation of a square and let $\tilde{S}$ be the fixed-corner function field defined above. For each triangle, let $\tilde{Z_i}\in \tilde{S}$ be the (possibly inhomogeneous) polynomial of degree less than or equal to $2$ expressing the area of this triangle in the $2k$ variables ${\tilde{X}_v, \tilde{Y}_v}$.   
Then there exists a homogeneous polynomial $f_T$ with integer coefficients in $n$ variables such that 
\begin{equation}\label{eq:fixcornf}
   2f_T(\tilde{Z_1},\ldots,\tilde{Z_n})=1 
\end{equation}
in $\tilde{S}$.  Furthermore, for $f_T$, we may take any polynomial satisfying the conclusion of Theorem Monsky+.  Conversely, any homogeneous polynomial of degree $e$ satisfying equation \ref{eq:fixcornf} also satisfies the conclusion of Theorem Monsky+.
\end{cor}

\begin{proof}
Let $f_T$ be a polynomial obtained from Theorem Monsky+.   For all vertices $v$, including corners, substitute $\tilde{X_v}, \tilde{Y}_v$ for $x_v, y_v$.  After these substitutions, each $Z_i$ becomes $\tilde{Z_i}$ and $\sigma$ becomes 1.  Hence equation \eqref{eq:fixcornf} is satisfied in $\tilde{S}$.

Conversely, suppose that $f_T$ is any homogeneous polynomial with integer coefficients satisfying Equation \eqref{eq:fixcornf}.  Let $M$ be the map defined in the proof of Theorem Monsky+, and for any vertex $v$, including corners, define $(\tilde{x_v}, \tilde{y_v})= M(x_v, y_v)$.   In equation \eqref{eq:fixcornf}, substituting $\tilde{x_v}, \tilde{y_v}$ for the variables ${\tilde{X}_v, \tilde{Y}_v}$ for all interior vertices $v$ turns each $\tilde{Z_i}$ corresponding to triangle $T_i=uvw$ into
\begin{equation*}
    \frac 12
\left|
\begin{matrix}
1 & 1 & 1 \\
\tilde{x}_u & \tilde{x}_v & \tilde{x}_w \\
\tilde{y}_u & \tilde{y}_v & \tilde{y}_w 
\end{matrix}
\right|
=
\frac1{2\sigma}
\left|
\begin{matrix}
1 & 1 & 1 \\
x_u & x_v & x_w \\
y_u & y_v & y_w 
\end{matrix}
\right|
= 
\frac1{\sigma} Z_i,
\end{equation*}
the first equation following from the fact that $M$ has determinant $1/\sigma$.  Hence with this substitution, equation \eqref{eq:fixcornf} becomes
$$
2f(\frac1{\sigma} Z_1, \dots, \frac1{\sigma} Z_n) = 1.
$$
Finally, by the homogeneity of $f$, we get the desired $2f( Z_1, \dots,  Z_n) = \sigma^e $. 
\end{proof}

\begin{example}[cf.~Example \ref{ex:2}]\label{ex:diag2}
Let $T$ be the triangulation shown in Figure \ref{fig:diag2}. Here there are 6 triangles with two interior vertices $u$ and $v$, so by Corollary \ref{cor:fixcorn} we may work in the rational function field in the 4 variables $x_u, y_u, x_v, y_v$.  The $Z$'s are defined by
\begin{align*}
    2Z_A &=  y_u    \\
    2Z_B &= x_v y_u-x_u y_v-y_u+y_v    \\
    2Z_C &= 1-x_v    \\
    2Z_D &=  x_u   \\
    2Z_E &=  x_u y_v-x_v y_u- x_u +x_v  \\
    2Z_F &= 1-y_v 
\end{align*} 
This time, there is no linear polynomial $f_T$ satisfying the desired conclusion.  However the quadratic
$$f_T = (A+C+E)^2-2AC+2DF +(A+C+E)(B+D+F)
$$ 
is a Monsky polynomial:  it has the property that $2f_T(Z_A, \dots, Z_F) = 1$.
\end{example}

\subsection{Constraints}\label{sec:++}

We next soup up our previous argument a bit more to take into account the constraints $\CC$.  
We assume $\T=(T,\CC)$ is drawable and that we have a parameterization $g$ of $X(\T)$ coming from a drawing order $\le$ as in Theorem \ref{param}.
We find that there is again a polynomial $f$ satisfying \eqref{eq:f}, this time with the areas of the living triangles expressed in terms of the parameters $w_i$ of the drawing order. 

Let $\T=(T,\CC)$ be a constrained triangulation of a square that is drawable. 
Let $\le$ be any drawing order, let $k=\sum\alpha_i$, and let $g_\le:\C^k\to X_\T$ be the parameterizing map defined earlier.
Denote the coordinates of $\C^k$ by $w_1,\ldots,w_k$.  Let $U=\Q(w_1, \dots, w_k)$ be the corresponding field of rational functions in $k$ variables.  We call $U$ the parameter field of the drawings of $T$.

Note that the number of living triangles is $n=k+2$.

\begin{thm}[Monsky++]
Let $\T=(T,\CC)$ be a constrained triangulation of a square that is drawable.   Fix a drawing order $\le$ with corresponding parameter field $U$. 
For each living triangle $\Delta_j$ ($1\le j\le n$), let $W_j\in U$ be the rational function in the $w_i$ expressing the area of $\Delta_j$, i.e., $W_j$ is the $j$th coordinate function of the map $\Area\circ\, g_\le$. Let $\sigma = \sum W_j$. 
Then there exists a homogeneous polynomial $f_\T$ with integer coefficients in $n$ variables such that
$$2f_\T(W_1,\ldots,W_n)=\sigma^e$$ 
in $U$, for some non-negative integer $e$.

In fact, if $f_T$ (note the font change) is the polynomial promised by Monsky+ for the honest triangulation $T$, then we may choose $f_\T$ to be the polynomial obtained from $f_T$ by plugging in zeroes for the variables that correspond to the dead triangles of $\T$.
\end{thm}

\begin{defn}
If $\T$ is a constrained triangulation of a square with $n$ living triangles, then any polynomial $f$ satisfying the conclusion of Theorem Monsky++ is called a \emph{Monsky polynomial} for $\T$.
\end{defn}

\begin{proof}
Fix $\T$ and $\le$.  
Let $m$ be the total number of triangles (alive and dead) in $T$.  
Recall $n=k+2=\sum \alpha_i +2$ is the number of living triangles.
We number the triangles of $T$ so that the first $k+2$ are living in $\T$.

Apply Monsky+ to the honest $T$ to get
$$
f_T(Z_1, \dots, Z_m) = \frac12 \sigma^e
$$
in the field $S= \Q(\{ x_i, y_i |\ 1\le i\le m\})$.
Here as above $Z_j$ is the polynomial in the $x_i, y_i$ expressing the area of $\Delta_j$.  
Let $f_\T$ denote the polynomial $f_T$ evaluated with all variables corresponding to dead triangles set to zero. 

We claim that
$$f_\T(W_1, \dots, W_{k+2}) = \frac12 (W_1+\dots+W_{k+2})^e$$
in $U$.  
To see this, specialize to any point $(w_i)$ in the domain of $g_{\leq}$, where the equation is an equality of complex numbers that is true by Monsky+, because these coordinates describe a drawing of $T$.
Since the domain of $g_{\leq}$ is dense in $\C^{(\sum \a)}$, the polynomials must be identical and the claim is established. 
\end{proof}

The preceding theorem could also be proved directly, in a manner very similar to the proof of Monsky+, but invoking the map $g$.  

These versions of Monsky's theorem provide the generalizations promised in the introduction.

\begin{cor}
If $\Diss$ is a generalized dissection of a parallelogram $\square$ into triangles with areas $a_1,\ldots,a_n$, then there is an integer polynomial $f$ in $n$ variables with $f(a_1,\ldots,a_n)=\frac 1 2 \Area(\square)^e$, for some non-negative integer $e$.  Moreover $f$ can be chosen to be invariant under deformation of $\Diss$.
\end{cor}

\begin{proof}
If the constrained triangulation $\T(\Diss)$ is drawable then we may apply Monsky++ directly to $\T(\Diss)$.  
In any case, even if $\T(\Diss)=(T,\CC)$ is not drawable, we apply Monsky+ to the honest triangulation $T$, getting the polynomial $f_T$.
The dissection $\Diss$ is the image of a drawing $\rho$; it doesn't matter that $\rho$ is not generic. 
As the areas of all dead triangles of $\T(\Diss)$ are zero in $\Diss$ and any deformation of $\Diss$, the polynomial $f_{\T(\Diss)}$ obtained from $f_T$ by plugging in zeroes for all variables corresponding to dead triangles in $\T(\Diss)$ satisfies the conclusion of the corollary.
\end{proof}

\begin{example}[$ACE$ again, cf.~Examples \ref{ex:2}, \ref{ex:2ACE-D}, \ref{ex:2ACE-T}]\label{ex:ACE++}
Let $\T$ be the $ACE$ example, i.e., the constrained triangulation shown in Figure \ref{fig:ACE-T}; a generic drawing is shown in Figure \ref{fig:ACE-D}.  Following the idea of Corollary \ref{cor:fixcorn}, we fix the corners to be the vertices of the unit square.  Here there are 3 living triangles called $A, C, E$ and two interior vertices $u$ and $v$.  We use the drawing order in which $u<v$ and see that $\alpha_u = 1$ while $\alpha_v=0$. Thus with fixed corners, there is just one parameter $w_1$.  The relevant condition for $u$ is $\SS\PP u$, and the relevant conditions for $v$ are $\RR\SS v$ and $\QQ vu$.  For any value of $w_1$, the drawing $g(w_1)$ places $u$ at $(0, 1-w_1)$ and places $v$ at $(\frac{w_1}{w_1-1} , 1)$.    The areas of the triangles are 
\begin{align*}
    W_A &= \frac12  (1-w_1)    \\
    W_C &=   \frac{1}{2(1-w_1)}  \\
    W_E &=   \frac{w_1^2}{2(w_1-1)}  
\end{align*} 
As the reader can see, $4W_A W_C = 1$, so the polynomial 
$$f_{\T}= 2AC$$ satisfies $2f_{\T}(W_A, W_C, W_E) =1$ and is a Monsky polynomial for $\T$.  Indeed, in agreement with Theorem Monsky++, this $f_{\T}$ equals the polynomial $f_T$ from Example \ref{ex:diag2} with the variables $B, D, F$ set to $0$.
\end{example}

\subsection{Computation}\label{sec:compute}

The last assertion of the Monsky++ theorem is that the polynomials $f_T$ and $f_\T$ are related in the most natural way possible.
(Of course these polynomials are not uniquely defined so this statement is not entirely precise.)
However in practice, to compute the polynomial $f_{\T(\Diss)}$ for a generalized dissection $\Diss$, we do not actually use this relationship.  
Doing that would require first invoking Monsky+, computing $f_T$ for the corresponding honest triangulation, and then zeroing out a bunch of variables.
But, as these are Gr\"obner basis computations which grow very quickly in complexity with the number of variables, it is far preferable to perform the calculation without introducing variables that we know we are eventually going to evaluate to zero.
This is the real value of the Monsky++ theorem:  it says that we can work directly with the parameters of the dissection, i.e., the coordinate functions of $g$, thereby reducing the variables to only those that are actually needed.
This is what we just saw in Example \ref{ex:ACE++}, where the deformation space has only one parameter.
As a result one can typically compute $f$ reasonably quickly for a generalized dissection with up to about 10 living triangles, even though the corresponding honest triangulation $T$ may have many more triangles than this and attempting to compute $f_T$ may crash our computers.
The left dissection of Figure \ref{fig:samples}, for example, has four parameters.  There are six triangles and its Monsky polynomial has degree four, so $f$ has at most $\binom{9}{4}= 126$ monomials and it is easily computed (in fact it has $104$ monomials).  The corresponding honest triangulation $T$ has 10 triangles and a Monsky polynomial of degree six.  This one is still computable in a reasonable amount of time, but it is quite large and unwieldy.

\section{The area variety}

Following \cite{triangles1}, we now introduce the machinery necessary to define the polynomial $p$, starting with the area map. 

Given three points $(x_1,y_1), (x_2,y_2), (x_3,y_3)\in\C^2$, we have defined the area of the oriented triangle $\Delta$ with these vertices (in this order) to be 
$$\Area(\Delta)=\frac 12
\left|
\begin{matrix}
1 & 1 & 1 \\
x_1 & x_2 & x_3 \\
y_1 & y_2 & y_3 
\end{matrix}
\right|.
$$
We also sometimes write $\Area(p_1p_2p_3)$ for the area of the triangle $\Delta$ with vertices $p_1,p_2,p_3\in\C^2$. 

Note that $\Area(p_1p_2p_3)=0$ if and only if $p_1,p_2,p_3$ lie on a (complex) line in $\C^2$.
When $\D\subset\R^2$ the function $\Area$ gives the usual (signed) area.

Let $\T$ be a fixed drawable constrained triangulation.  
Let $T_1,\dots,T_n$ be the living triangles of $\T$. 
Note each $T_i$ inherits an orientation from $T$.
Let $Y(\T)$ be the projective space $\P^{n-1}$ with coordinates $[\cdots:A_i:\cdots]$, $1\le i\le n$.
If $\T=\T(\Diss)$ then we also denote $Y(\T)$ by $Y(\Diss)$.

As $\T$ is drawable, the set of generic drawings is open and dense in $X(\T)$. 
We therefore have a rational map 
$$\Area=\Area_\T:X(\T)\dashrightarrow Y(\T)$$
given by
$\Area(\rho)=[\cdots:\Area(T_i):\cdots]$.

\begin{defn}
Let $\T$ be drawable.  The \emph{area variety} $V=V(\T)$ is the closure in $Y(\T)$ of $\Area(X(\T))$.
If $\T=\T(\Diss)$ for some generalized dissection $\Diss$ then we also refer to $V=V(\T)=V(\Diss)$ as the area variety of $\Diss$.
\end{defn}

Note that the affine group $\Aff_2(\C)$ acts on $X$ and that the area map is equivariant with respect to this action.  
(This accounts for not just the translations and scaling we alluded to after defining $X$, but also rotations and shears.)  
Thus since $\T$ is drawable the generic fibers of the area map are at least 6-dimensional (that being the size of $\Aff_2$).

Meanwhile one easily counts that $\dim(Y)=n-1=\dim X-5$.

Therefore, for a given $\T$ the area variety $V(\T)$ has codimension at least 1 in $Y(\T)$, with equality if and only if generic fibers are exactly 6-dimensional.  
In other words, $V$ is a hypersurface in $Y$ if and only if there is no 1-parameter family of area-preserving deformations, other than those contained in an $\Aff_2$ orbit.   

\begin{defn}
A drawable constrained triangulation $\T$ is \emph{hyper} if $V(\T)$ is a hypersurface in $Y(\T)$.
\end{defn}

\begin{conjecture}\label{conj:punt2}
If $\Diss$ is a generalized dissection then $\T(\Diss)$ is hyper, i.e., $V(\Diss)$ is a hypersurface in $Y(\Diss)$.
\end{conjecture}

At least two phenomena can prevent an arbitrary constrained triangulation from being hyper, as we described in Section 4 of \cite{triangles1}.\footnote{In each of those scenarios there is a subset of the variables that sums to zero and as far as we know the area variety is still a hypersurface in a smaller projective space than $Y$.  So, it is possible that a slight modification of Conjecture \ref{conj:punt2} holds for all $\T$.}
However these phenomena do not arise for constrained triangulations of the form $\T(\Diss)$.

All honest triangulations are hyper, and we proved in \cite{triangles1} that if $\T$ has only one constraint and it is non-separating then $\T$ is hyper.
That proof can be extended somewhat.  In a forthcoming paper we further enlarge the set of $\T(\Diss)$'s that we know to be hyper. 

If $\T$ is hyper then there is a unique (up to scaling)
non-zero polynomial $p=p_\T$ that vanishes on $V$.  
The polynomial $p$ is irreducible because $X$, and
therefore $V$, is an irreducible variety. 
Also $p$ has rational coefficients (because the coordinate functions of $\Area$ do) so $p$ can be normalized to have integer coordinates with no common factor.  
We assume this has been done; the polynomial $p_\T$ is now well-defined up to sign for any hyper $\T$.

\begin{defn}
We call $p$ and $-p$ the \emph{area polynomials} for $\T$.
\end{defn}

We remark that the area polynomials are computable, using Gr\"obner basis techniques, but that these computations quickly become intractable as the triangulation grows.

\section{Mod 2}\label{sec:mod2}

Let $\T$ be a constrained triangulation that is hyper (hence drawable).
We  thus have Monsky polynomials $f$ and an area polynomial $p$, both homogeneous elements of the polynomial ring $\Z[A_1, \dots, A_n]$, with variables $A_i$ corresponding to the living triangles of $\T$.  Letting $\sigma= \sum A_i$, the equations $2f=\sigma^e$ and $p=0$ hold on the variety $V(\T)$.  The existence of $p$ and $f$ with these properties is enough to reveal all of the coefficients of $p$ modulo 2.

\begin{thm}[Mod 2 theorem]\label{thm:mod2}
If $\T$ is hyper then its area polynomial $p$ satisfies $$p \equiv \sigma^d \mod 2,$$
where $d=\deg p$.
\end{thm}

\begin{proof}
Since $V(\T)$ is the zero set of the irreducible polynomial $p$, any polynomial that vanishes on $V(\T)$ is a multiple of $p$.  Thus $p\ |\ 2f-\sigma^e$.  
These are integer polynomials and the divisibility occurs in $\Q[A_1, \dots, A_n]$, so there is a polynomial $q$ in this ring with
$p\cdot q=2f-\sigma^e$.  
By Gauss's Lemma, $q$ must have integer coefficients.  
Therefore we may reduce the coefficients mod 2, and using $[\cdot]$ for the reduction, we have $[p][q]= [\sigma^e]$.
Since $\Z/2\Z[A_1,\ldots,A_n]$ is a unique factorization domain, we conclude that $[p]=[\sigma^d]$ for some $d=\deg p\le e$.
\end{proof}

\begin{cor}
Let $\T$ be hyper, and suppose the area polynomial $p$ has degree $d$.  Then all leading terms $A_i^d$ occur with odd (hence non-zero) coefficient. 
\end{cor}

In a forthcoming paper we show that the leading coefficients are all equal up to sign.  We suspect these coefficients are all $\pm 1$, as we discuss in Section \ref{sec:positivity}.

\section{Canonical Monsky polynomials}
\label{sec:canonical}

Let $\T$ be hyper, let $p$ be its area polynomial, and suppose $p$ has degree $d$.
Recall that $p$ is only well-defined up to sign; we suppose we have chosen one of the two possibilities.

Using $p$, we now single out a particular $f$ satisfying Monsky++.  By Theorem \ref{thm:mod2} we have 
$\sigma^d + p =2f$
for some polynomial $f\in\Z[A_1, ..., A_n]$, homogeneous of degree $d=\deg p$.
Note that $f$ satisfies the conclusion of Monsky++, since $2f-\sigma^d=p$ vanishes on $V$.
This shows that we may choose $$f=\frac 1 2 (\sigma^d+p)$$ in Theorem Monsky++.  This shows in addition that we may choose a Monsky polynomial with the same degree as $p$ (and no lower).  

If we begin with $-p$ instead of $p$, we end up with 
$$\tilde f=\frac 1 2 (\sigma^d-p)=f-p$$ 
instead.  The pair $\{f,\tilde f\}$ is therefore a canonically defined pair of Monsky polynomials, both of minimal degree.  

\begin{defn}
For $\T$ hyper with area polynomial $\pm p$ of degree $d$, the \emph{canonical Monsky polynomials} are $\{f,\tilde f\}$ where $f=\frac 1 2 (\s^d+p)$ and $\tilde f=\frac 1 2 (\s^d-p)$.
\end{defn}

\begin{prop}
For $\T$ hyper with area polynomials $\pm p$ of degree $d$, the polynomial $f_0$ is a Monsky polynomial for $\T$ if and only if $f_0=\frac 1 2 (\s^e+pq)$ for some integer $e\ge d$ and for some polynomial $q\equiv \sigma^{e-d}\mod 2$.

In particular, the minimal degree Monsky polynomials are exactly $\{f+mp\}=\{\ftilde+mp\}$, where $f,\ftilde$ are the canonical Monsky polynomials and where $m$ is any integer. 
\end{prop}

\begin{proof}
For polynomials $f_0\in\Z[A_1, ..., A_n]$, the condition that $f_0$ be a Monsky polynomial is equivalent to the condition that $p\ |\ 2f_0 - \sigma^e$ for some $e$, which in turn is equivalent to the condition that $f_0=\frac 1 2 (\sigma^e +pq)$ for some $q\equiv \sigma^{e-d}\mod 2$.  The second assertion follows by considering $e=d$.
\end{proof}

The canonical Monsky polynomials $f,\tilde f$ satisfy 
\begin{align}\label{eq:fp}
    f + \ftilde & = \s^d
\\  f - \ftilde & = p
\end{align}
(where interchanging $f$ and $\tilde f$ corresponds to interchanging $p$ and $-p$).

\begin{example}[cf.~Examples \ref{ex:1}, \ref{ex:1-D}, \ref{ex:diag1}]
We saw  in Example \ref{ex:diag1} that for the triangulation in Figure \ref{fig:diag1}, $A+C$ is a Monsky polynomial. 
The area polynomial is $p=A-B+C-D$ and so the canonical Monsky polynomials are $f=B+D$ and $\tilde{f}=A+C$.  We have $2f=2\tilde f= \s$ in the field $U$ and on the variety $V$.  Other Monsky polynomials may be obtained by adding a multiple of $p$ to $f$.  For example, $p+f=2A-B+2C-D$ is a Monsky polynomial. 
\end{example}

\begin{example}[cf.~Example \ref{ex:2}, \ref{ex:diag2}]
For the honest triangulation $T$ of Figure \ref{fig:diag2} we choose the area polynomial $p=(A+C+E)^2-4AC-(B+D+F)^2+4DF$.  (See also \cite{triangles1}, where this is worked out in detail.)  This is irreducible and it vanishes on $V$.  The canonical Monsky polynomials are 
\begin{align*}
    f &=\frac 1 2 (\s^2+p)=(A+C+E)^2-2AC+2DF + (A+C+E)(B+D+F)\\
    \ftilde &=\frac 1 2 (\s^2-p)= (B+D+F)^2-2DF+2AC +(A+C+E)(B+D+F).
\end{align*}
This is how we found the polynomial $f_T$ given in Example \ref{ex:2} from the introduction and again in Example \ref{ex:diag2} in Section \ref{sec:+}.

Notice in this example that there is another way to obtain $f$ and $\tilde f$ from $p$.  If we write $p=p_+-p_-$ as the difference of two polynomials with non-negative coefficients and no common terms, then we have $p_+=p_-$ on $V$, where
\begin{align*}
p_+&=A^2+C^2+E^2+2AE+2CE+2DF \qquad \mbox{ and}\\ p_-&=B^2+D^2+F^2+2BD+2BF+2AC.    
\end{align*}
We see that each of these coefficients is less than or equal to the corresponding coefficient in the expansion of $\s^2$.  The terms of this expansion that do not occur in $p_+$ or $p_-$ are $2(A+C+E)(B+D+F)$.  Thus we can ``make up the difference'' by adding $t=(A+C+E)(B+D+F)$ to both sides, giving $p_++t=p_-+t$ on $V$, and 
\begin{align}\label{eq:t}
    (p_++t) + (p_-+t) & = \s^d
\\  (p_++t) - (p_-+t) & = p.
\end{align}
This means we have found two polynomials whose sum is $\s^d$ and whose difference vanishes on $V$; therefore each is half of $\s^d$ and they are Monsky polynomials.
Comparing with \eqref{eq:fp}, we see that in fact 
\begin{align*}
    f &= p_++t \qquad \mbox{ and}\\
    \ftilde &=p_-+t.
\end{align*}
\end{example}

\section{Positivity}\label{sec:positivity}

Something happened in the last derivation that may not always work.  When we wrote $p_+$ and $p_-$, we observed that all the terms have coefficients that are ``small,'' in the sense that none is larger than the corresponding coefficient of $\s^2$.  As a direct consequence, all coefficients of $t$, hence also of both $f$ and $\tilde f$, turned out to be non-negative.

This phenomenon is not essential to the procedure; if it fails one can still define $t$ satisfying \eqref{eq:t} and proceed to determine $f$ and $\tilde f$.  In that case $t$ and at least one of $f,\tilde f$ would have some negative coefficients.
For honest triangulations, however, we have never observed this.  

\begin{defn}
A polynomial is called \emph{positive} if all its coefficients are non-negative.
\end{defn}

\begin{conjecture}[Positivity]
The canonical Monsky polynomials of every honest triangulation are positive.
\end{conjecture}

We can restate this conjecture in terms of the area polynomial as follows.  

\begin{defn}
A homogeneous polynomial $p$ of degree $d$ is called \emph{\good} if each coefficient of $p$ has absolute value less than or equal to the corresponding multinomial coefficient; that is, if both $\s^d-p$ and $\s^d+p$ are positive.
\end{defn}

Because of the relationships $2f=\s^d+p$ and $2\tilde f=\s^d-p$, the positivity conjecture is equivalent to saying that the area polynomial of an honest triangulation is \good.

At the end of the previous section we mentioned our suspicion that the leading terms $A_i^d$ of the area polynomial $p$ all have coefficient $\pm 1$.  We point out now that this is a special case of the positivity conjecture. 

We conclude this section with some further remarks about this conjecture.
We frame these remarks in terms of the area polynomial $p$, rather than $f$, because we have more techniques for computing and working with $p$.

Observe that \goodness\ is preserved under products:  if $p$ and $\bar p$ are two \good\ polynomials, then $p \bar p$ is also \good.

Observe also that \goodness\ is a local condition on the polynomial $p$, in the sense that its failure is always witnessed by (at least) one individual monomial.
In other words $p$ is \good\ if and only if each monomial of $p$ is \good, even though a given monomial may not include all the variables in the polynomial $p$.  For instance any polynomial containing the term, say, $-40ABCD$ fails to be \good, because the degree is 4 and $|\!-\!40|$ is larger than the coefficient of $ABCD$ in $\s^4$, which is $24$ regardless of how many variables there are in $\s$.

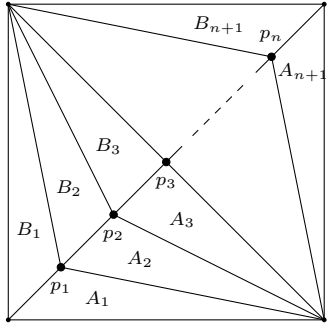
\begin{figure}
    \begin{tikzpicture}[scale=.7]
            \draw (0,6) rectangle (6,0) ;
\draw[fill=black] (0,0) circle (1pt);
\draw[fill=black] (0,6) circle (1pt);
\draw[fill=black] (6,6) circle (1pt);
\draw[fill=black] (6,0) circle (1pt);
\draw (0,0) -- (3,3);
\draw (3,3) -- (3.3,3.3);
\draw[dashed] (3.5,3.5) -- (4.5,4.5);
\draw (4.7,4.7) -- (5,5);
\draw (5,5) -- (6,6);
\draw (0,6) -- (1,1);
\draw (0,6) -- (2,2);
\draw (0,6) -- (3,3);
\draw (0,6) -- (5,5);
\draw (6,0) -- (1,1);
\draw (6,0) -- (2,2);
\draw (6,0) -- (3,3);
\draw (6,0) -- (5,5);
\draw[fill=black] (1,1) circle (2pt);
\draw[fill=black] (2,2) circle (2pt);
\draw[fill=black] (3,3) circle (2pt);
\draw[fill=black] (5,5) circle (2pt);
\begin{scriptsize}
\draw (1.7,0.4) node{$A_1$};
\draw (2.5,1.15) node{$A_2$};
\draw (3.3,1.9) node{$A_3$};
\draw (5.6,4.7) node{$A_{n+1}$};
\draw (0.4,1.7) node{$B_1$};
\draw (1.15,2.5) node{$B_2$};
\draw (1.9,3.3) node{$B_3$};
\draw (4,5.6) node{$B_{n+1}$};
\draw (1,1) node[label=below:$p_1$] {};
\draw (2,2) node[label=below:$p_2$] {};
\draw (3,3) node[label=below:$p_3$] {};
\draw (5,5) node[label=above:$p_n$] {};
\end{scriptsize}
    \end{tikzpicture}
    \caption{The diagonal case $\T_n$.}
    \label{fig:diagonalcase}
\end{figure}

With these observations, and using the methods of \cite{triangles1}, we can show that the positivity conjecture holds for the infinite family $\T_n$ of honest triangulations shown in Figure \ref{fig:diagonalcase}.  We call this the ``diagonal case.''  

Note that $\T_1$ and $\T_2$ have already made numerous appearances in this paper under the pseudonyms Example 1 and Example 2.

\begin{prop} 
For each $n$, the diagonal case $\T_n$ has positive Monsky polynomial.
\end{prop}

\begin{proof}
As it is honest, $\T_n$ is hyper, and we denote its area polynomial by $p_n$.  In \cite{triangles1} we gave an explicit expression for $p_n$, but it is difficult to tell directly from that expression that $p_n$ is \good.  However we do know that the degree of $p_n$ is exactly $n$.  
If we focus on any particular monomial of $p_n$, then at least one subscript from $\{1,\ldots,n+1\}$, call it $j$, does not occur in the variables of this monomial. If we kill triangles $A_j$ and $B_j$ the resulting polynomial $p_n|_{A_j=B_j=0}$ factors into a linear factor with coefficients $\pm 1$ and a degree $n-1$ factor which is a version of $p_{n-1}$ (but with subscripts at least $j$ shifted up by one).  The linear factor is \good, and inductively so is $p_{n-1}$, so their product is also small.  The monomial from $p_n$ on which we focused is one term of this product, and so this chosen monomial satisfies \goodness.  Since our choice of monomial was arbitrary, it follows that $p_n$ is \good, as desired.
\end{proof}

\begin{example}[$ACE$ again, cf.~Examples \ref{ex:2}, \ref{ex:2ACE-D}, \ref{ex:2ACE-T}, \ref{ex:ACE++}]
To see some non-honest examples, one can start with the honest triangulation $\T_2$ and plug in zeroes.  For instance consider the $ACE$ example.  Monsky++ implies that we may take the $f$ and $\tilde f$ computed already for the honest case and plug in $B=D=F=0$, giving 
\begin{align*}
    f &= A^2+C^2+E^2+2AE+2CE\\
    \ftilde &=2AC.
\end{align*}
The area polynomial $p=(A+C+E)^2-4AC$ is likewise obtained by plugging in $B=D=F=0$ from Example \ref{ex:diag2}.
We check $f+\ftilde=(A+C+E)^2$ and $f-\ftilde=p$.
\end{example}

There is a reason that we assume honesty in the positivity conjecture.  The next example shows a classical dissection $D$ whose $\T(D)$ has an area polynomial that is not \good.

\begin{example}[Failure of positivity]
We return to the second figure in this paper, reproduced in Figure \ref{fig:nonpositive-D}.  The constrained triangulation $\T=\T(D)=(T,\CC)$ has ten triangles and four constraints; there are six living triangles.  The area polynomial $p_T$ for the honest triangulation $T$ has degree six, and is \good.  However, the area polynomial $p_\T$ for the constrained triangulation has degree four and has a total of 70 terms, one of which is $-40ABCD$, where $A,B,C,D$ denote the areas of the four triangles that touch the corners.  Therefore $p_\T$ is not \good.  Of course, this term does not occur in the degree six polynomial $p_T$.

Curiously, of the 70 terms of $p_\T$,  only one violates \goodness.
Likewise, of the 104 terms of $f$ and the 122 terms of $\tilde f$, just one of them has a negative coefficient, namely $-8ABCD$.

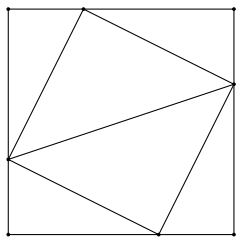
\begin{figure}
    \begin{center} 
    \begin{tikzpicture}[scale=.5]
        \draw (0,0) -- (0,6) -- (6,6) -- (6,0) -- cycle;
\draw (0,2) -- (2,6) -- (6,4) -- (4,0) -- cycle;
\draw (0,2) -- (6,4) ;
\draw[fill=black] (0,0) circle (1pt);
\draw[fill=black] (0,6) circle (1pt);
\draw[fill=black] (6,0) circle (1pt);
\draw[fill=black] (6,6) circle (1pt);
\draw[fill=black] (0,2) circle (1pt);
\draw[fill=black] (2,6) circle (1pt);
\draw[fill=black] (6,4) circle (1pt);
\draw[fill=black] (4,0) circle (1pt);
    \end{tikzpicture} 
    \caption{A $\T$ whose Monsky polynomial is not positive.}
    \label{fig:nonpositive-D}
\end{center}
\end{figure}
\end{example}

\section{Equidissections}\label{sec:equi}

\begin{figure}
    \centering
    \begin{tikzpicture}[scale=.4]
            \draw (0,0) -- (0,6) -- (6,6) -- (6,0) -- cycle;
\draw (3,1.6) -- (4,4) -- (1.6,3.6) -- cycle;
\draw (4,4) -- (6,0) -- (3,1.6) -- (0,0) -- (1.6,3.6) -- (0,6) -- (4,4) -- (6,6) ;
\draw[fill=black] (0,0) circle (1pt);
\draw[fill=black] (0,6) circle (1pt);
\draw[fill=black] (6,0) circle (1pt);
\draw[fill=black] (6,6) circle (1pt);
\draw[fill=black] (3,1.6) circle (1pt);
\draw[fill=black] (4,4) circle (1pt);
\draw[fill=black] (1.6,3.6) circle (1pt);
    \end{tikzpicture}
    \caption{This dissection cannot be deformed to an equidissection.}
    \label{fig:3int}
\end{figure}
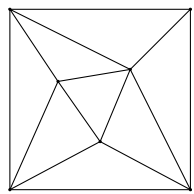

Monsky's original equidissection theorem leaves open the question of which dissections can be deformed to be equi-areal.  (It is easy to see that for each even $n$, there exist such dissections with $n$ triangles.) 
Observe that if $D$ is a dissection with $n$ triangles and $p(1,1,\ldots,1)\ne0$, or equivalently the sum of the coefficients of $f$ not equal to $n^d/2$, then $D$ cannot be deformed to an equidissection without killing triangles.
This is the case, for instance, with the dissection shown in Figure \ref{fig:3int}:  there are $8$ triangles, and the polynomial $p$ has degree $3$, so plugging in all $1$'s to $\s^d$ gives the value $8^3=512$.  However plugging in all $1$'s in $f$ and $\tilde f$ gives the values 260 and 252, which are not equal (and $p(1,\ldots,1)$ is the difference, $\pm 8$).
There is no drawing of this triangulation in which all triangles have area $1/8$.  In this case, this is also easily proved using elementary Euclidean geometry. 

\begin{question}
Given a dissection or generalized dissection $\Diss$, can one predict from the combinatorics of $\T(\Diss)$ whether or not $\Diss$ can be deformed to an equidissection?
\end{question}

\section*{acknowledgements}
 The authors thank the Budapest Semesters in Mathematics program for supporting us in 2019 as Director's Mathematicians in Residence.  The first author acknowledges additional support from the Simons Foundation (grant \#281189) and the Lenfest Summer Research Fund of Washington and Lee University.
 
We are also grateful to Paul Monsky for sharing his valuable insights and perspectives.

\bibliographystyle{plain}
\bibliography{TriangleRefs}


\end{document}